\documentclass[11pt]{article}

\usepackage{amsmath,amssymb,amsfonts}
\usepackage{pstricks,pst-node}

\numberwithin{equation}{section}
\newtheorem{thm}{Theorem}[section] 
\newtheorem{prp}[thm]{Proposition}
\newtheorem{lmm}[thm]{Lemma}   
\newtheorem{crl}[thm]{Corollary}

\def\e_ref#1{(\ref{#1})}

\def\ov#1{\overline{#1}}
\def\ti#1{\tilde{#1}}

\def\smsize#1{\begin{small}#1\end{small}}

\def\lra{\longrightarrow}

\def\eset{\emptyset}
\def\i{\infty}

\def\io{\iota}

\def\om{\omega}
\def\si{\sigma}
\def\vp{\varphi}

\def\C{\mathcal C}

\def\D{\mathcal D}

\def\R{\mathbb R}
\def\Z{\mathbb Z}

\def\De{\Delta}

\def\Bd{\textnormal{Bd}~\!}

\def\CH{\textnormal{CH}}
\def\ev{\textnormal{ev}}
\def\id{\textnormal{id}}
\def\Id{\textnormal{Id}}
\def\Im{\textnormal{Im}~\!}
\def\Int{\textnormal{Int}~\!}

\def\rk{\textnormal{rk}~\!}
\def\sd{\textnormal{sd}\,}
\def\sign{\textnormal{sign}~\!}
\def\top{\textnormal{top}}

\def\cal{\mathcal}
\def\Bbb{\mathbb}

\begin{document}

\title{Pseudocycles and Integral Homology} 
\author{Aleksey Zinger}
\date{\today}
\maketitle

\begin{abstract}
\noindent
We describe a natural isomorphism between 
the set of equivalence classes of pseudocycles and 
the integral homology groups of a smooth manifold.
Our arguments generalize to settings well-suited for applications 
in enumerative algebraic geometry and for construction of 
the virtual fundamental class in the Gromov-Witten theory.
\end{abstract}

\thispagestyle{empty}
\tableofcontents

\section{Introduction}

\subsection{Main Theorem}

\noindent
In his seminal paper~\cite{G}, 
Gromov initiated the study of pseudoholomorphic curves 
in symplectic manifolds and demonstrated their usefulness by proving 
a number of important results in symplectic topology.
In~\cite{McSa} and~\cite{RT}, pseudoholomorphic curves are used to 
define invariants of semipositive manifolds.
In particular, it is shown in~\cite{McSa} and~\cite{RT}
that for every compact semipositive symplectic manifold~$(X,\om)$,
homology class $A\!\in\!H_2(X;\Bbb{Z})$,
integers $k\!\ge\!3$ and $N\!\ge\!0$, and 
generic compatible almost complex structure $J$ on~$X$,
there exists a smooth oriented manifold  ${\cal M}_{k,N}(A,J)$
and a smooth map
$$\ev_{k,N}^{A,J}\!: {\cal M}_{k,N}(A,J)\lra X^{k+N}$$
such that the ``boundary'' of $\ev_{k,N}^{A,J}$ is small; see below.
Such a smooth map is called a {\it pseudocycle}
and determines a homomorphism $H_*(X^{k+N};\Bbb{Z})\!\lra\!\Bbb{Z}$,
which turns out to be an invariant of~$(X,\om;A,k,N)$.\\

\noindent
In general, if $X$ is a smooth manifold, subset $Z$ of $X$ 
is said to have {\it dimension at most~k} 
if there exists a \hbox{$k$-dimensional} manifold~$Y$ and 
a smooth map $h\!:Y\!\lra\!X$ such that the image of~$h$ contains~$Z$.
If $f\!:M\!\lra\!X$ is a continuous map between topological spaces,
the {\it boundary} of~$f$ is the~set
$$\Bd f=\bigcap_{K\subset M\text{~cmpt}}\!\!\!\!\!\!\ov{f(M\!-\!K)}.$$
A smooth map $f\!:M\!\lra\!X$ is a {\it $k$-pseudocycle} if
$M$ is an oriented $k$-manifold, $f(M)$ is a 
pre-compact\footnote{i.e.~its closure is compact} subset of~$X$, and
the dimension of $\Bd f$ is at \hbox{most $k\!-\!2$}.
Two $k$-pseudocycles $f_0\!:M_0\!\lra\!X$ and $f_1\!:M_1\!\lra\!X$
are {\it equivalent} if there exists a smooth oriented manifold~$\ti{M}$
and a smooth map $\ti{f}\!:\ti{M}\!\lra\!X$ such that
the image of~$\ti{f}$ is a pre-compact subset of~$X$,
$$\dim\Bd\ti{f}\le k\!-\!1,\quad \partial\ti{M}=M_1-M_0,
\quad \ti{f}|_{M_0}=f_0, \quad\hbox{and}\quad \ti{f}|_{M_1}=f_1.$$
We denote the set of equivalence classes of pseudocycles into $X$
by~${\cal H}_*(X)$.
This set is naturally a $\Z$-graded module over~$\Z$.
In this paper, we prove

\begin{thm}
\label{main_thm}
If $X$ is a smooth manifold, there exist 
natural\,\footnote{In other words,  $\Psi_*$ and $\Phi_*$ are natural transformations
of functors $H_*(\cdot;\Z)$ and ${\cal H}_*(\cdot)$ from
the category of smooth compact manifolds and maps.}
homomorphisms of graded $\Z$-modules
$$\Psi_*\!: H_*(X;\Z)\lra{\cal H}_*(X)
\quad\hbox{and}\quad
\Phi_*\!: {\cal H}_*(X)\lra H_*(X;\Z),$$
such that $\Phi_*\circ\Psi_*\!=\!\Id$ and $\Psi_*\circ\Phi_*\!=\!\Id$.
\end{thm}

\noindent
{\it Remark 1:} In~\cite{McSa} and~\cite{RT}, a pseudocycle is not explicitly
required to have a pre-compact image.
As~\cite{McSa} and~\cite{RT} work with compact manifolds,
this condition is automatically satisfied.
However, this requirement is essential in the non-compact case.
As observed in~\cite{K}, there is no surjective homomorphism from
$H_*(X;\Z)$ to ${\cal H}_*(X)$ if $X$ is not compact 
and pseudocycles are not required to have pre-compact images.\\

\noindent
{\it Remark 2:} It is sufficient to require that a pseudocycle map be continuous,
as long as the same condition is imposed on pseudocycle equivalences.
All arguments in this paper go through for continuous pseudocycles.
In fact, Lemma~\ref{bdpush_lmm} is no longer necessary. 
However, smooth pseudocycles are useful in symplectic topology 
for defining invariants as intersection numbers.\\

\noindent
In order to define symplectic invariants, 
\cite{McSa} and~\cite{RT} observe that every element of~${\cal H}_*(X)$ 
defines a homomorphism $H_*(X;\Bbb{Z})\!\lra\!\Bbb{Z}$, or
equivalently an element of 
$H_*(X;\Bbb{Z})\big/\hbox{Tor}\big(H_*(X;\Bbb{Z})\big)$.
Thus, Theorem~\ref{main_thm} leads to symplectic invariants
that may be strictly stronger than the GW-invariants defined 
in \cite{McSa} and in~\cite{RT}.
In fact, they are as good as the maps~$\ev_{k,N}^{A,J}$ can give:

\begin{crl}
If $(X,\om_1)$ and $(X,\om_2)$ are semipositive
symplectic manifolds that have the same \hbox{$GW$-invariants},
viewed as a collection of integral homology classes,
then the corresponding collections of evaluation maps 
from products of moduli spaces of pseudoholomorphic maps 
and Riemannian surfaces are equivalent as pseudocycles.
\end{crl}

\noindent
This corollary is immediate from Theorem~\ref{main_thm}.\\

\noindent
This paper was begun while the author was a graduate student at MIT
and then put on the back burner.
Its aim was to clarify relations between $H_*(X;\Z)$ and~${\cal H}_*(X)$
that were hinted at in~\cite{McSa} and stated without a proof in~\cite{RT}.
Since then, this issue has been explored in~\cite{K} and in~\cite{Sc}.
The views taken in  \cite{K} and in~\cite{Sc} differ significantly from the present paper.
In particular, while non-compact manifolds are considered in~\cite{K},
pseudocycles in~\cite{K} are not required to have pre-compact preimages.
Theorem~\ref{main_thm} fails for such pseudocycles.
The arguments in the present paper are rather direct and use no advanced techniques,
beyond standard algebraic topology.
In a sense they implement an outline proposed in Section~7.1 of~\cite{McSa}.
However, a fully rigorous implementation of this outline requires non-trivial technical 
facts obtained in Subsections~\ref{neigh_subs} and~\ref{singular} of this paper.
An additional obstacle arises in showing that the map $\Psi_*$ is well-defined.
As discussed in the next subsection, there are two ways of overcoming it.\\

\noindent
As a graduate student, the author was partially supported by an 
NSF Graduate Research Fellowship and NSF grant DMS-9803166.
The author would like to thank D.~McDuff for a suggestion that greatly simplified 
one of the steps in the original argument and for comments on an earlier draft
of the present version that led to improvements in the exposition.

\subsection{Outline of Constructions}

\noindent
If $M$ is a compact oriented $k$-manifold and $f\!:M\!\lra\!X$ is a smooth map,
$f$ determines an element of $H_k(X;\Z)$, i.e.~the pushforward of the fundamental
class of $M$, $f_*[M]$.
If $f\!:M\!\lra\!X$ is a $k$-pseudocycle, $M$ need not be compact.
However, one can choose a compact $k$-submanifold with boundary, $\bar{V}\!\subset\!M$,
so that $f(M\!-\!V)$ lies in a small neighborhood $U$ of~$\Bd f$.
In particular, $f|_{\bar{V}}$ determines the homology class
$$f_*[\bar{V}]\in H_k(X,U;\Z).$$
By Proposition~\ref{neighb_prp}, $U$ can be chosen so that $H_k(X,U;\Z)$ is
naturally isomorphic to $H_k(X;\Z)$.
In order to show that the image of $f_*[V]$ in $H_k(X;\Z)$ depends only on $f$
(and not $V$ or $U$), we use  Proposition~\ref{isom1_prp}  to replace
the singular chain complex~$S_*(X)$  by a quotient complex~$\bar{S}_*(X)$.
The advantage of the latter complex is that cycles and boundaries between chains
can be constructed more easily.\\

\noindent
In an analogous way, a pseudocycle equivalence $\ti{f}\!:\ti{M}\!\lra\!X$ between
two pseudocycles 
$$f_i\!:M_i\!\lra\!X, \qquad i\!=\!0,1,$$ 
gives rise to a chain equivalence 
between the corresponding cycles in $\bar{S}_*(X,W)$, for a small neighborhood $W$
of $\Bd\ti{f}$.
In particular, the homology cycles determined by $f_0$ and $f_1$ are 
equal in $H_k(X,W;\Z)$.
On the other hand, by Proposition~\ref{neighb_prp}, $W$ can be chosen so that 
$H_k(X;\Z)$ naturally injects into $H_k(X,W;\Z)$.
Therefore, $f_0$ and $f_1$ determine the same elements of $H_k(X;\Z)$ and
the homomorphism $\Phi_*$ is well-defined.
Its construction is described in detail in Subsection~\ref{phi_sec}.\\

\noindent
{\it Remark 1:} The homomorphism 
$$\Phi_*\!:{\cal H}_*(X)\lra H_*(X;\Z)$$
of Subsection~\ref{phi_sec} induces the linear map
$${\cal H}_*(X) \lra H_*(X;\Z)
\big/\hbox{Tor}\big(H_*(X;\Z)\big)$$
described in~\cite{McSa} and~\cite{RT}.
However, our construction of $\Phi_*$ differs from that of the induced map
in~\cite{McSa} and~\cite{RT}.
Indeed, the latter is constructed via the homomorphism $\Psi_*$
and a natural intersection pairing on~${\cal H}_*(X)$.
The construction of~$\Phi_*$ in Subsection~\ref{phi_sec}
is more direct.\\

\noindent
{\it Remark 2:} The construction of Subsection~\ref{phi_sec} implies the following.
Suppose $M$ is an oriented $k$-manifold and $f\!:M\!\lra\!X$ is a continuous map 
with a pre-compact image. 
Suppose $\Bd f$ has an arbitrary small neighborhood~$U$ so that 
$$H_l(U)=0 \qquad\forall~l> k\!-\!2.$$
Then, $f$ defines an element in 
$H_k(X)$\footnote{This statement holds for any coefficient ring.}.
An analogous statement holds for equivalences between maps.
Note that we are not assuming that $X$ is a smooth manifold.
These observations have a variety of applications.
For example, the first statement implies that a compact complex algebraic
variety carries a fundamental class.
For essentially the same reason, (generalized) pseudocycles figure 
prominently in the approach in~\cite{genus0pr} to a large class of
problems in enumerative geometry.
Pseudocycles can also be used to give a more geometric interpretation
of the virtual fundamental class construction of~\cite{FO} and~\cite{LT}
and are used to define new symplectic invariants in~\cite{g1comp2}.
This is a different type of generalization, as
the ambient space $X$ in these settings is a topological space stratified
by infinite-dimensional manifolds.\\

\noindent
In order to construct the homomorphism $\Psi_*$, we would like to show that 
a singular cycle gives rise to a pseudocycle and a chain equivalence between 
two cycles gives rise to a pseudocycle equivalence between the corresponding
pseudocycles.
The former works out precisely as outlined in Section~7.1 of~\cite{McSa},
with a reinterpretation for the complex 
$\bar{S}_*(X)$\footnote{This reinterpretation is not necessary to construct 
the map $\Psi_*$ in Subsection~\ref{psi_sec}.
However, it is needed in  Subsection~\ref{isom_sec} to show that the maps
$\Phi_*$ and $\Psi_*$ are isomorphisms.}.
If $s$ is a $k$-cycle, all codimension-one simplices of its $k$-simplices
must cancel in pairs.
By gluing the $k$-simplices along the codimension-one faces paired in this way,
we obtain a continuous map from a compact topological space $M'$ to~$X$.
The complement of the codimension-two simplices is a smooth manifold
and the continuous map can be smoothed out in a fixed manner
using Lemma~\ref{bdpush_lmm}.
We thus obtain a pseudocycle from the cycle~$s$.\\

\noindent
On the other hand, turning a chain equivalence $\ti{s}$ between two $k$-cycles,
$s_0$ and $s_1$, 
into a pseudocycle equivalence between the corresponding pseudocycles,
$$f_0\!:M_0\lra X \qquad\hbox{and}\qquad f_1\!:M_1\lra X,$$
turns out to be less straightforward.
Similarly to the previous paragraph, $\ti{s}$ gives rise to a smooth map
from a smooth $(k\!+\!1)$-manifold with boundary,
$$\ti{f}\!:\ti{M}^*\lra X.$$
However, if all codimension-two simplices (including those of dimension $k\!-\!1$)
are dropped, the boundary of $\ti{M}^*$ will be the complement in 
$M_0\!\sqcup\!M_1$ of a subset of dimension~$k\!-\!1$
(instead of being $M_0\!\sqcup\!M_1$).
One way to fix this is to keep the $(k\!-\!1)$-simplices that would lie 
on the boundary.
In such a case, the entire space may no longer be a smooth manifold and 
its boundary may not be $M_0\!\sqcup\!M_1$, because the $(k\!-\!1)$-simplices
of the $(k\!+\!1)$-simplices of $\ti{s}$ may be identified differently
from the way  the $(k\!-\!1)$-simplices of the $k$-simplices of $s_0$
and $s_1$ are identified.
It is possible to modify $\ti{s}$ so that all identifications are consistent. 
However, the required modification turns out to be quite laborious.
This direct approach was implemented in the original version of this paper 
(now called {\it Version~B}).\\

\noindent
In the present version, we instead implement a less direct, but far simpler, construction
which is based on a suggestion of D.~McDuff.
Instead of trying to reinsert $(k\!-\!1)$-simplices into the boundary of $\ti{M}^*$,
we will attach to $\ti{M}^*$ two collars, 
$$\ti{M}_0\subset [0,1]\!\times\!M_0  \qquad\hbox{and}\qquad
 \ti{M}_1\subset [0,1]\!\times\!M_1.$$
The boundary of $\ti{M}_i$ will have two pieces, $M_i$ and the complement in 
$M_i$ of the $(k\!-\!1)$-simplices.
We attach the latter to the piece of the boundary of $\ti{M}^*$ corresponding
to~$s_i$.
In this way, we obtain a smooth manifold $\ti{M}$ with boundary $M_1\!-\!M_0$
and a pseudocycle equivalence from $f_0$ to~$f_1$; 
see Subsection~\ref{psi_sec} for details.\\

\noindent
Finally, in Subsection~\ref{isom_sec}, we verify that the homomorphisms $\Psi_*$
and $\Phi_*$ are mutual inverses.
It is fairly straightforward to see that the map $\Phi_*\!\circ\!\Psi_*$ is
the identity on $\bar{H}_*(X;\Z)$.
However, showing the injectivity of $\Phi_*$ requires more care.
The desired pseudocycle equivalence,
$$\ti{f}\!:\ti{M}\lra X,$$
is constructed by taking a limit of the corresponding construction 
in Subsection~\ref{psi_sec}.
In particular, the smooth manifold $\ti{M}$ is obtained  as a subspace of
a {\it non}-compact space.

\section{Preliminaries}

\subsection{Notation}

\noindent
In this subsection we describe a number of subsets of 
the standard $k$-simplex $\De^k$ as well as maps between standard simplices
of various dimensions.
We conclude by constructing self-maps of $\De^k$ which will be used to 
define canonical smoothings of piecewise smooth maps in Subsection~\ref{psi_sec}; 
see Lemma~\ref{bdpush_lmm}.
The reader may want to skip this subsection at first and refer back for notation
as needed.\\

\noindent
If $A$ is a finite subset of $\R^k$, we denote by $\CH(A)$ and $\CH^0(A)$ 
the (closed) convex hull of~$A$ and the open convex hull of~$A$, 
respectively, i.e.
\begin{equation*}\begin{split}
CH(A)&=\Big\{\sum_{v\in A}t_vv\!: t_v\!\in\![0,1];\,\sum_{v\in A}t_v\!=\!1\Big\}
\qquad\hbox{and}\\
CH^0(A)&=\Big\{\sum_{v\in A}t_vv\!: t_v\!\in\!(0,1);\,\sum_{v\in A}t_v\!=\!1\Big\}.
\end{split}\end{equation*}\\
For each $p\!=\!1,\ldots,k$, let $e_p$ be the $p$th coordinate vector in $\R^k$.
Put $e_0\!=\!0\!\in\!\R^k$.
Denote by
$$\De^k=\CH\big(e_0,e_1,\ldots,e_k\big) \qquad\hbox{and}\qquad
\Int\De^k=\CH^0\big(e_0,e_1,\ldots,e_k\big)$$
the standard $k$-simplex and its interior.
Let 
$$b_k=\frac{1}{k\!+\!1}\bigg(\sum_{q=0}^{q=k}e_q\bigg)=
\Big(\frac{1}{k\!+\!1},\ldots,\frac{1}{k\!+\!1}\Big)\in\Bbb{R}^k.$$
be the barycenter of $\De^k$.\\

\noindent
For each $p\!=\!0,1,\ldots,k$, let
$$\De^k_p=\CH\big(\big\{e_q\!: q\!\in\!\{0,1,\ldots,k\}\!-\!p\big\}\big)
\qquad\hbox{and}\qquad
\Int\De_p^k=\CH^0\big(\big\{e_q\!: q\!\in\!\{0,1,\ldots,k\}\!-\!p\big\}\big)$$
denote the $p$th face of $\De^k$ and its interior.
Define a linear map\footnote{A map $f\!:\De^m\!\lra\!\De^k$ is called {\it linear} if
$$f(t_0e_0\!+\!\ldots\!+t_me_m)= t_0f(e_0)\!+\!\ldots\!+\!t_mf(e_m)
\qquad\forall~(t_0,\ldots,t_m)\in[0,1]^{m+1}~~\hbox{s.t.}~~
t_0\!+\!\ldots\!+\!t_m\!=\!1.$$}
$$\io_{k,p}\!: \De^{k-1}\lra\De^k_p\subset\De^k \qquad\hbox{by}\qquad
\io_{k,p}(e_q)=\begin{cases}
e_q,& \hbox{if}~q\!<\!p;\\
e_{q+1},& \hbox{if}~q\!\ge\!p.
\end{cases}$$
We also define a projection map
$$\ti\pi^k_p\!: \De^k\!-\!\{e_p\}\lra\De^k_p \qquad\hbox{by}\qquad
\ti\pi^k_p\Big(\sum_{q=0}^{q=k}t_qe_q\Big)=
\frac{1}{1\!-\!t_p}\Big(\sum_{q\neq p}t_qe_q\Big).$$
Put
$$b_{k,p}=\io_{k,p}(b_{k-1}),\qquad
b_{k,p}'=\frac{1}{k\!+\!1}\bigg(b_k+\!\sum_{q\neq p}e_q\bigg).$$
The points $b_{k,p}$ and $b_{k,p}'$ are the barycenters of 
the $(k\!-\!1)$-simplex~$\De^k_p$
and of the $k$-simplex spanned by~$b_k$ and the vertices of~$\De^k_p$.
Define a neighborhood of $\Int\De^k_p$ in $\De^k$ by
$$U^k_p=\big\{t_pb_{k,p}'\!+\!\sum_{0\le q\le k;q\neq p}\!\!\!\!\!\!\!t_qe_q\!: 
t_p\!\ge\!0,~t_q\!>\!0~\forall q\!\neq\!p;~ \sum_{q=0}^{q=k}t_q\!=\!1\big\}.$$
These neighborhoods will be used to construct pseudocycles out of homology cycles.\\

\begin{figure}
\begin{pspicture}(-.7,-1.8)(10,1.5)
\psset{unit=.4cm}
\psline{->}(0,-3)(8,-3)\psline{->}(1,-4)(1,4)\psline(1,3)(7,-3)
\pscircle*(1,-3){.2}\pscircle*(7,-3){.2}\pscircle*(1,3){.2}
\rput(.4,-3.5){\smsize{$e_0$}}\rput(7.4,-3.5){\smsize{$e_1$}}
\rput(.3,2.7){\smsize{$e_2$}}
\rput(4,-3.7){\smsize{$\De^2_2$}}\rput{90}(.3,0){\smsize{$\De^2_1$}}
\rput{-45}(5,0){\smsize{$\De^2_0$}}
\pscircle*(3,-1){.2}\pscircle*(3.67,-2.33){.2}
\rput(2.4,-1){\smsize{$b_2$}}\rput(4.6,-1.9){\smsize{$b_{2,2}'$}}
\psline{->}(15,-3)(23,-3)\psline{->}(16,-4)(16,4)\psline(16,3)(22,-3)
\pscircle*(16,-3){.2}\pscircle*(22,-3){.2}\pscircle*(16,3){.2}
\rput(15.4,-3.5){\smsize{$e_0$}}\rput(22.4,-3.5){\smsize{$e_1$}}
\rput(15.3,2.7){\smsize{$e_2$}}\rput(19,-3.7){\smsize{$\De^2_2$}}
\pscircle*(16,0){.2}\pscircle*(19,0){.2}
\rput(15.1,-.3){\smsize{$b_{2,1}$}}\rput(20.2,0){\smsize{$b_{2,0}$}}
\psline{->}(30,-3)(38,-3)\psline{->}(31,-4)(31,4)\psline(31,3)(37,-3)
\pscircle*(31,-3){.2}\pscircle*(37,-3){.2}\pscircle*(31,3){.2}
\rput(30.4,-3.5){\smsize{$e_0$}}\rput(37.4,-3.5){\smsize{$e_1$}}
\rput(30.3,2.7){\smsize{$e_2$}}
\rput(34,-3.7){\smsize{$\De^2_2$}}
\psline[linestyle=dotted](31,-3)(33.67,-2.33)\psline[linestyle=dotted](37,-3)(33.67,-2.33)
\pnode(33.5,-2.65){A}\rput(29,0){\rnode{B}{\smsize{$U_2^2$}}}
\nccurve[nodesep=0,angleA=-90,angleB=180,ncurv=.5,linewidth=.02]{->}{B}{A}
\end{pspicture}
\caption{The Standard $2$-Simplex and Some of its Distinguished Subsets}
\label{simplices_fig}
\end{figure}

\noindent
If $p,q\!=\!0,1,\ldots,k$ and $p\!\neq\!q$, let 
$$\De^k_{p,q}\equiv \De^k_p\!\cap\!\De^k_q$$
be the corresponding codimension-two simplex.
Define neighborhoods of $\Int\De^k_{p,q}$ in $\De^k$ by
\begin{equation*}\begin{split}
\ti{U}^k_{p,q}&=\big\{t_pb_{k,p}\!+\!t_qb_{k,q}\!+\!\sum_{r\neq p,q}\!\!t_re_r\!:  
t_p,t_q\!\ge\!0,~t_r\!>\!0~\forall\,r\!\neq\!p,q;~ \sum_{r=0}^{r=k}t_r\!=\!1\big\},\\  
U^k_{p,q}&=\big\{t_p\io_{k,p}(b_{k-1,\io_{k,\io_{k,p}^{-1}(q)}}')
\!+\!t_q\io_{k,q}(b_{k-1,\io_{k,\io_{k,q}^{-1}(p)}}')
\!+\!\sum_{r\neq p,q}\!\!t_re_r\!:
t_p,t_q\!\ge\!0,\, t_r\!>\!0\, \forall r\!\neq\!p,q;\,\sum_{r=0}^{r=k}t_r\!=\!1\big\}.
\end{split}\end{equation*}
These sets will be used to construct pseudocycle equivalences out of 
singular chains.
Note that
\begin{equation}\label{disjoint_e}
\ti{U}^k_{p_1,q_1}\cap\ti{U}^k_{p_2,q_2}=\eset  \qquad\hbox{if}\quad
\{p_1,q_1\}\neq\{p_2,q_2\}.
\end{equation}
Define a projection map
$$\ti\pi^k_{p,q}\!: \De^k-\CH(e_p,e_q)\lra\De^k_{p,q} \qquad\hbox{by}\qquad
\ti\pi^k_{p,q}\Big(\sum_{r=0}^{r=k}t_re_r\Big)=
\frac{1}{1\!-\!t_p\!-\!t_q}\Big(\sum_{r\neq p,q}t_re_r\Big).$$\\

\begin{figure}
\begin{pspicture}(-8.7,-1.8)(10,1.5)
\psset{unit=.4cm}
\psline{->}(-12,-3)(-4,-3)\psline{->}(-11,-4)(-11,4)\psline(-11,3)(-5,-3)
\pscircle*(-11,-3){.2}\pscircle*(-5,-3){.2}\pscircle*(-11,3){.2}
\rput(-11.6,-3.5){\smsize{$e_0$}}\rput(-4.6,-3.5){\smsize{$e_1$}}
\rput(-11.7,2.7){\smsize{$e_2$}}
\pscircle*(-8,-3){.2}\pscircle*(-11,0){.2}
\psline[linestyle=dotted](-8,-3)(-11,0)\rput(-10.1,-2.2){\smsize{$\ti{U}^2_{2,1}$}}
\psline{->}(2,-3)(10,-3)\psline{->}(3,-4)(3,4)\psline(3,3)(9,-3)
\pscircle*(3,-3){.2}\pscircle*(9,-3){.2}\pscircle*(3,3){.2}
\rput(2.4,-3.5){\smsize{$e_0$}}\rput(9.4,-3.5){\smsize{$e_1$}}
\rput(2.3,2.7){\smsize{$e_2$}}
\psline[linestyle=dotted](3,-1.5)(4.5,-3)
\pnode(3.5,-2.5){A}\rput(1,0){\rnode{B}{\smsize{$U^2_{2,1}$}}}
\nccurve[nodesep=0,angleA=-90,angleB=180,ncurv=.5,linewidth=.02]{->}{B}{A}
\end{pspicture}
\caption{The Standard $2$-Simplex and Some of its Distinguished Subsets}
\label{simplices_fig2}
\end{figure}

\noindent
Finally, let ${\cal S}_k$ denote the group of permutations of the set
$\{0,\ldots,k\}$. 
The set ${\cal S}_k$ can be viewed as a subset of ${\cal S}_{k+1}$; 
if $\tau\!\in\!{\cal S}_k$, put \hbox{$\tau(k\!+\!1)\!=\!k\!+\!1$}. 
For any $\tau\!\in\!{\cal S}_k$, let
$$\tau\!:\De^k\lra\De^k$$
be the linear map defined by
$$\tau(e_q)=e_{\tau(q)} \qquad\forall\,q=0,1,\ldots,k.$$

\begin{lmm}
\label{bdpush_lmm}
If $k\!\ge\!1$, $Y$ is the $(k\!-\!2)$-skeleton of $\De^k$, and 
$\ti{Y}$ is the $(k\!-\!2)$-skeleton of $\De^{k+1}$,
there exist continuous functions 
$$\vp_k\!: \De^k\lra\De^k \qquad\hbox{and}\qquad 
\ti\vp_{k+1}\!: \De^{k+1}\lra\De^{k+1}$$
such that\\
${}~~~~$ (i) $\vp_k$ is smooth outside of $Y$ and 
$\ti\vp_{k+1}$ is smooth outside of $\ti{Y}$;\\
${}~~~$ (ii) for all $p\!=\!0,\ldots,k$ and $\tau\!\in\!{\cal S}_k$,
\begin{equation}\label{bdpush_e1}
\vp_k|_{U^k_p}=\ti\pi^k_p\big|_{U^k_p} \qquad\hbox{and}\qquad
\vp_k\circ\tau=\tau\circ\vp_k;
\end{equation}
${}~~$ (iii) for all $p,q\!=\!0,\ldots,k\!+\!1$ with $p\!\neq\!q$ and 
$\ti\tau\!\in\!{\cal S}_{k+1}$,
\begin{equation}\label{bdpush_e2}
\ti\vp_{k+1}|_{U^{k+1}_{p,q}}=\ti\pi^{k+1}_{p,q}\big|_{U^{k+1}_{p,q}}, \qquad
\ti\vp_{k+1}\circ\ti\tau=\ti\tau\circ\ti\vp_{k+1}, \quad\hbox{and}\quad
\ti\vp_{k+1}\circ\io_{k+1,p}=\io_{k+1,p}\circ\vp_k.
\end{equation}\\
\end{lmm}

\noindent
{\it Proof:} (1) Choose a smooth function
$$\ti\eta_{0,1}\!: \De^{k+1}-\De^{k+1}_{0,1}\!\cap\!\ti{Y}\lra [0,1]$$
such that  $\ti\eta_{0,1}\!=\!1$ on $U^{k+1}_{0,1}$, 
$\ti\eta_{0,1}\!=\!0$ outside of~$\ti{U}^{k+1}_{0,1}$, and
$\ti\eta_{0,1}$ is invariant under any permutation $\ti\tau\!\in\!{\cal S}_{k+1}$
that preserves the set~$\{0,1\}$.
If $\ti\tau\!\in\!{\cal S}_{k+1}$ is any permutation, let
$$\ti\eta_{\ti\tau(0),\ti\tau(1)}=\ti\eta_{0,1}\!\circ\!\ti\tau^{-1}\!:
\De^{k+1}-\De^{k+1}_{\ti\tau(0),\ti\tau(1)}\!\cap\!\ti{Y}\lra [0,1].$$
By the assumptions on $\ti\eta_{0,1}$, $\ti\eta_{p,q}$ is a well-defined smooth
function such that $\ti\eta_{p,q}\!=\!1$ on $U^{k+1}_{p,q}$, 
$\ti\eta_{p,q}\!=\!0$ outside of~$\ti{U}^{k+1}_{p,q}$, and
\begin{equation}\label{etainv_e}
\ti\eta_{\ti\tau(p),\ti\tau(q)}=\ti\eta_{p,q}\circ\ti\tau^{-1}
\end{equation}
for all $\ti\tau\!\in\!{\cal S}_{k+1}$ and distinct $p,q\!=\!0,\ldots,k\!+\!1$.\\

\noindent
(2) Define
$$\ti\vp_{k+1}\!:\De^{k+1}\lra\De^{k+1} \qquad\hbox{by}\qquad
\ti\vp_{k+1}(x)=x+\!\!\!\sum_{0\le p<q\le k+1}\!\!\!\!\!\!
\ti\eta_{p,q}(x)\cdot\big(\ti\pi^{k+1}_{p,q}(x)\!-\!x\big).$$
Since $\ti\eta_{p,q}$ vanishes on a neighborhood of $\CH(e_p,e_q)$
and $\ti\pi_{p,q}^{k+1}$ restricts to the identity on~$\De^{k+1}_{p,q}$, 
the function~$\ti\eta$ is well-defined,
continuous everywhere, and smooth on $\De^{k+1}\!-\!\ti{Y}$.
By~\e_ref{disjoint_e}, $\ti\vp_{k+1}\!=\!\ti\pi^{k+1}_{p,q}$ on~$U^{k+1}_{p,q}$.
By~\e_ref{etainv_e}, for every $\ti\tau\!\in\!{\cal S}_{k+1}$
\begin{equation*}\begin{split}
\ti\vp_{k+1}\circ\ti\tau
&=\ti\tau+\!\!\!\sum_{0\le p<q\le k+1}\!\!\!\!\!\!
\ti\eta_{p,q}\!\circ\!\ti\tau\cdot\big(\ti\pi^{k+1}_{p,q}\!\circ\!\ti\tau\!-\!\ti\tau\big)\\
&=\ti\tau+\!\!\!\sum_{0\le p<q\le k+1}\!\!\!\!\!\!
\ti\eta_{\ti\tau^{-1}(p),\ti\tau^{-1}(q)}\cdot
\big(\ti\tau\!\circ\!\ti\pi^{k+1}_{\ti\tau^{-1}(p),\ti\tau^{-1}(q)}\!-\!\ti\tau\big)\\
&=\ti\tau+\!\!\!\sum_{0\le p<q\le k+1}\!\!\!\!\!\!
\ti\eta_{p,q}\cdot\big(\ti\tau\!\circ\!\ti\pi^{k+1}_{p,q}\!-\!\ti\tau\big)
=\ti\tau\circ\ti\vp_{k+1}.
\end{split}\end{equation*}
Thus, $\ti\vp_{k+1}$ satisfies the first two conditions in~\e_ref{bdpush_e2}, 
as well as~(i) above.\\

\noindent
(3)  We define~$\vp_k$ by the third condition in~\e_ref{bdpush_e2}.
The function~$\vp_k$ is independent of the choice of~$p$ and satisfies 
the second condition in~\e_ref{bdpush_e1}.
To see this, suppose $p,q\!=\!0,\ldots,k\!+\!1$ and $\tau\!\in\!{\cal S}_k$.
Let $\ti\tau\!\in\!{\cal S}_{k+1}$ be defined~by 
$$\ti\tau\circ\io_{k+1,p}=\io_{k+1,q}\circ\tau.$$
If $\vp_{k,p}$ and $\vp_{k,q}$ are the functions corresponding to~$p$ 
and~$q$ via the third equation in~\e_ref{bdpush_e2}, then by 
the second equation in~\e_ref{bdpush_e2}
\begin{equation*}\begin{split}
\io_{k+1,q}\circ\tau\circ\vp_{k,p}
=\ti\tau\circ\io_{k+1,p}\circ\vp_{k,p}
&=\ti\tau\circ\ti\vp_{k+1}\circ\io_{k+1,p}
=\ti\vp_{k+1}\circ\ti\tau\circ\io_{k+1,p}\\
&=\ti\vp_{k+1}\circ\io_{k+1,q}\circ\tau
=\io_{k+1,q}\circ\vp_{k,q}\circ\tau.
\end{split}\end{equation*}
We conclude that
$$\tau\circ\vp_{k,p}=\vp_{k,q}\circ\tau
\qquad\forall~p,q\!=\!0,\ldots,k\!+\!1,~\tau\!\in\!{\cal S}_k.$$
The function $\vp_k$ satisfies the first condition in~\e_ref{bdpush_e1} because
\begin{gather*}
\io_{k+1,p}\big(U^k_p\big)=U^{k+1}_{p,p+1}\cap\De^{k+1}_{p,p+1} \qquad\hbox{and}\\
\io_{k+1,p}\circ\vp_k=\ti\vp_{k+1}\circ\io_{k+1,p}
=\ti\pi^{k+1}_{p,p+1}\circ\io_{k+1,p}=\io_{k+1,p}\circ\ti\pi^k_p
~~~\hbox{on}~~U^k_p.
\end{gather*} 
Finally, $\vp_k$ satisfies (i) because $\ti\vp_{k+1}$ does.

\subsection{Homology of Neighborhoods of Smooth Maps}
\label{neigh_subs}

\noindent 
In this subsection, we prove

\begin{prp}
\label{neighb_prp}
If $h\!: Y\!\lra\!X$ is a smooth map and $W$ is an open neighborhood of
a subset $A$ of $\Im h$ in $X$, 
there exists a neighborhood $U$ of $A$ in $W$ such that
$$H_l(U)=0 \qquad\hbox{if}~~l\!>\!\dim Y.$$\\
\end{prp}

\noindent
Note that it may not be true that
$$H_l(A)=0 \qquad \hbox{if}~~l\!>\!\dim Y.$$
For example, let $A$ be the subset of $X\!=\!\Bbb{R}^N$ 
consisting of countably many $k$-spheres of radii tending to~$0$ and 
having a single point in common.
If $k\!\ge\!2$, the set~$A$ has infinitely many nonzero
homology groups, as shown in~\cite{BM}.\\

\noindent
If $h\!: Y\!\lra\!X$ is a smooth map and
$k$ is a nonnegative integer, put
$$N_k(h)=\big\{y\!\in\!Y\!: \rk dh|_y\!\le\!k\big\}.$$
Proposition~\ref{neighb_prp} follows from Lemma~\ref{neighb_lmm}
applied with $X$ replaced by $W$, $Y$ by $h^{-1}(W)$,
and $k$ \hbox{by $\dim Y$}.\\

\noindent
One of the ingredients in the proof of Lemma~\ref{neighb_lmm}
is Lemma~\ref{transverse_lmm}.
For the purposes of this paper,
a {\it triangulation} of a smooth manifold~$X$ is a pair
$T\!=\!(K,\eta)$ consisting of a simplicial complex
and a homeomorphism $\eta\!:|K|\!\lra\!X$,
where $|K|$ is a geometric realization of~$K$ in~$\Bbb{R}^N$
in the sense of Section~3 in~\cite{Mu2}, 
such that $\eta|_{\Int\si}$ is smooth for every \hbox{simplex $\si\!\in\!K$}.

\begin{lmm}
\label{transverse_lmm}
If $X,Y$ are smooth manifolds and $h\!:Y\!\lra\!X$ is a smooth map,
there exists a triangulation $T\!=\!(K,\eta)$ of $X$ 
such that $h$ is transverse to $\eta|_{\Int\si}$
for every \hbox{simplex $\si\!\in\!K$}.
\end{lmm}

\noindent
This lemma is clear. 
In fact, we can start with any triangulation of~$X$ and 
obtain a desired one by an arbitrary small generic perturbation.

\begin{lmm}
\label{neighb_lmm}
If $h\!: Y\!\lra\!X$ is a smooth map, for every nonnegative integer~$k$
there exists a neighborhood $U$ of $h\big(N_k(h)\big)$ in $X$ such that
$$H_l(U)=0 \qquad\hbox{if}~~l\!>\!k.$$
\end{lmm}

\noindent
{\it Proof:}
By Lemma~\ref{transverse_lmm}, there exists a triangulation 
$T\!=\!(K,\eta)$ of $X$ such that the smooth map~$h$ is transversal to 
$\eta|_{\Int\si}$  for all $\si\!\in\!K$. 
In particular,
$$h\big(N_k(h)\big)\subset\!
\bigcup_{\si\in K,\dim\si\ge n-k}\!\!\!\!\!\!\!\!\!\!\!\!\!\!\eta(\Int\si)
=\!\bigcup_{\si\in K,\dim\si\ge n-k}\!\!\!\!\!\!\!\!\!\!\!\!\!\!
\eta\big(\hbox{St}(b_{\si},\sd K)\big),$$
where $n\!=\!\dim X$, $\sd K$ is the barycentric subdivision of~$K$,
and $\hbox{St}(b_{\si},\sd K)$ is the star of~$b_{\si}$
in~$\sd K$.\footnote{If $K$ is a simplicial complex and $\si$ is a simplex in $K$,
{\it the star of $\si$ in $K$} is the union of the subsets $\Int\,\si'$ taken
over the simplices $\si'\!\in\!K$ such that $\si\!\subset\!\si'$;
see Section~62 in~\cite{Mu2}.}
Note~that
$$\hbox{St}(b_{\si},\sd K)\cap\hbox{St}(b_{\si'},\sd K)=\eset$$
unless $\si\!\subset\!\si'$ or $\si'\!\subset\!\si$.
Furthermore, if $\si_1\!\subset\!\ldots\!\subset\!\si_m$,
$$\hbox{St}(b_{\si_1},\sd K)\cap\ldots\cap
\hbox{St}(b_{\si_m},\sd K)=
\hbox{St}(b_{\si_1}\ldots b_{\si_m},\sd K);$$
the last set is contractable. Put 
$$U_m'=\!\bigcup_{\si\in K,\dim\si=m}\!\!\!\!\!\!\!\!
\!\!\hbox{St}(b_{\si},\sd K).$$
Then $U_{l_m}'\!\cap\!\ldots\!\cap\!U_{m_j}'$
is a disjoint union of contractable open sets in~$|K|$.
Let 
$$U_m=\eta(U_m'), \quad m=n\!-\!k,\ldots,n;
\qquad
U=\bigcup_{m=n-k}^n\!\!\!U_m.$$
Since $\eta\!: |K|\!\lra\!X$ is a homeomorphism, 
$U_{m_1}\!\cap\!\ldots\!\cap\!U_{m_j}$
is a disjoint union of contractable open subsets of~$X$.
It follows  from Lemma~\ref{mv} below that $H_l(U)\!=\!0$ if $l\!>\!k$.
Furthermore, by the above \hbox{$h\big(N_k(h)\big)\!\subset\!U$}.

\begin{lmm}
\label{mv}
Let $\big\{U_m\big\}_{m=0}^{m=k}$ be a collection of open sets in $X$
and $U\!=\!\bigcup\limits_{m=0}\limits^{m=k}U_m$. 
If 
$$H_l\big(U_{m_1}\!\cap\!\ldots\!\cap\!U_{m_j};\Z\big)=0
\qquad\forall~l\!>\!0,~m_1,\ldots,m_j\!=\!0,\ldots,k,$$  
then $H_l(U)\!=\!0$  $l\!>\!k$.
\end{lmm}

\noindent
This lemma follows by induction from Mayer-Vietoris Theorem;
see~\cite[p186]{Mu2}.

\subsection{Oriented Homology Groups}
\label{singular}

\noindent
If $X$ is a simplicial complex, the standard singular chain complex $S_*(X)$
most naturally corresponds to the {\it ordered} simplicial chain complex of~$X$;
see Section~13 in~\cite{Mu2}.
In this subsection, we define a singular chain complex~$\bar{S}_*(X)$
which corresponds to the standard, or {\it oriented}, simplicial chain complex.
In particular, its homology is the same as the homology of 
the ordinary singular chain complex; see Proposition~\ref{isom1_prp}.
On the other hand, it is much easier to construct cycles
in~$\bar{S}_*(X)$ than in~$S_*(X)$.\\

\noindent
If $X$ is a topological space, let $\big(S_*(X),\partial_X\big)$ denote 
its singular chain complex, i.e.~the free abelian group
on the~set 
$$\bigcup_{k=0}^{\i}\hbox{Hom}(\De^k,X)$$
of all continuous maps from standard simplices to $X$,
along with a map $\partial_X$ of degree $-1$.
Let $S_k'(X)$ denote the free subgroup of $S_*(X)$ spanned by the set
$$\big\{f\!-\!(\sign\tau)f\!\circ\!\tau\!: 
f\!\in\!\hbox{Hom}(\De^k,X);~
\tau\!\in\!{\cal S}_k;~ k\!=\!0,1,\ldots\big\}.$$
If $\tau\!\in\!{\cal S}_k$, put
$$\ti\tau=\Id_{\De^k}-(\sign\tau)\tau\in S_k'(\De^k).$$
Then, $S_*'(X)$ is the subgroup of $S_*(X)$ spanned by 
$$\big\{f_{\#}\ti\tau\!: f\!\in\!\hbox{Hom}(\De^k,X);~
\tau\!\in\!{\cal S}_k;~ k\!=\!0,1,\ldots\big\}.$$
Note that if $h\!: X\!\lra\!Y$ is a continuous map, the linear~map 
$$h_{\#}\!: S_*(X)\lra S_*(Y)$$
 maps $S_*'(X)$ into~$S_*'(Y)$.

\begin{lmm}
\label{complex}
The free abelian group $S_*'(X)$ is a subcomplex of 
$\big(S_*(X),\partial_X\big)$,
i.e. $\partial_XS_*'(X)\!\subset\!S_*'(X)$.
\end{lmm}

\noindent
{\it Proof:} Suppose $\tau\!\in\!{\cal S}_k$. 
For any $p\!=\!0,\ldots,k$, let $\tau_p\!\in\!{\cal S}_{k-1}$ be such that
$$\tau\circ\io_{k,p}\!=\!\io_{k,\tau(p)}\circ\tau_p\!:
\De^{k-1}\lra\De^k_p\!\subset\!\De^k.$$
Let $\tau_{k,p}\!\in\!{\cal S}_k$ be defined by
$$\tau_{k,p}(q)
=\begin{cases}
\io_{k,p}(q),& \hbox{if}~q\!<\!k;\\
p,& \hbox{if}~q\!=\!k,\\
\end{cases}$$
Then,
$\tau\!\circ\!\tau_{k,p}=\tau_{k,\tau(p)}\!\circ\!\tau_p\in{\cal S}_k$
for all $\tau\!\in\!{\cal S}_k$.
Thus,
\begin{equation}\label{signch_e}
\sign\tau_p=(-1)^{(k-p)+(k-\tau(p))}\sign~\tau
=(-1)^{p+\tau(p)}\sign\tau.
\end{equation}
(2) By the above, we have
\begin{equation*}\begin{split}
\partial_{\De^k}\tau
&=\sum_{p=0}^k(-1)^p\tau\!\circ\!\io_{k,p}
=\sum_{p=0}^k(-1)^p\io_{k,\tau(p)}\!\circ\!\tau_p\\
&=(\sign\tau)\sum_{p=0}^k(-1)^{\tau(p)}(\sign\tau_p)
\io_{k,\tau(p)}\!\circ\!\tau_p
=(\sign\tau)\sum_{p=0}^k(-1)^p(\sign\tau_{\tau^{-1}(p)})
\io_{k,p}\!\circ\!\tau_{\tau^{-1}(p)}.
\end{split}\end{equation*}
Thus,
$$\partial_{\De^k}\tilde{\tau}=
\sum_{p=0}^k(-1)^p\big(\io_{k,p}-(\sign\tau_{\tau^{-1}(p)})
\io_{k,p}\!\circ\!\tau_{\tau^{-1}(p)}\big)\in S_{k-1}'(\De^k).$$
It follows that for any $f\!\in\!S_k(X)$, 
$$\partial_X(f_{\#}\ti\tau)=
f_{\#}(\partial_{\De^k}\ti\tau)\in S_{k-1}'(X).$$

\begin{lmm}
\label{homot}
There exists a natural transformation of functors
$D_X\!: S_*\!\lra\! S_{*+1}$ such that\\
(i) if $f\!: \De^m\!\lra\!\De^k$ is a linear map, 
$D_Xf$ is a linear combination of linear maps 
$\De^{m+1}\!\lra\!\De^k$ for all $k,m\!=\!0,1,\ldots$;\\
(ii) $D_XS_*'(X)\subset S_*'(X)$ for all topological spaces $X$;\\
(iii) $\partial_XD_X=(-1)^{k+1}\Id+D_X\partial_X$ on $S_k'(X)$.
\end{lmm}

\noindent
{\it Proof:}
(1) Suppose $k\!\in\Z^+$.
If $f\!:\De^m\!\lra\!\De^k$ is a linear map, 
define a new linear~map
\begin{equation}\label{homot_e1}
P_kf\!:\De^{m+1}\lra\De^k  \qquad\hbox{by}\qquad
P_kf\,(e_q)=
\begin{cases}
f(e_q),& \hbox{if}~q\!=\!0,\ldots,m;\\
b_k,& \hbox{if}~q\!=\!m\!+\!1.
\end{cases}
\end{equation}
The transformation~$P_k$ induces a linear map on the subchain complex of 
$S_*(\De^k)$ spanned by the linear maps.
If $\tau\!\in\!{\cal S}_m\!\subset\!{\cal S}_{m+1}$ and $f\!\in\!S_m(\De^k)$, 
then
\begin{equation}\label{homot_e1b}
P_k(f\!\circ\!\tau)=P_kf\circ\tau.
\end{equation}
Thus, $P_k$ maps the subgroup of $S_*'(\De^k)$ spanned by the linear maps 
into itself.
Similarly, if $\tau\!\in\!{\cal S}_k$, 
\begin{equation}\label{homot_e2}
\tau_{\#}(P_kf)\equiv\tau\circ P_kf=P_k(\tau\!\circ\!f)
\equiv P_k(\tau_{\#}f).
\end{equation}
Furthermore,
\begin{equation}\label{homot_e3}
\partial_{\De^k}P_kf=(-1)^{m+1}f+P_k\big(\partial_{\De^k}f\big).
\end{equation}\\

\noindent
(2) Let $D_X|_{S_k(X)}\!=\!0$ if $k\!<\!1$;
then $D_X$ satisfies (i)-(iii).
Suppose $k\!\ge\!1$ and we have defined $D_X|_{S_{k-1}(X)}$ so 
that the three requirements are satisfied wherever $D_X$ is defined.
Put
\begin{equation}\label{homot_e4}
D_{\De^k}(\Id_{\De^k})=
P_k\big(\Id_{\De^k}+(-1)^{k+1}D_{\De^k}\partial_{\De^k} \Id_{\De^k}\big)
\in S_{k+1}(\De^k).
\end{equation}
By the inductive assumption~(i) and \e_ref{homot_e1},
$D_{\De^k}(\Id_{\De^k})$ is a well-defined linear combination of linear maps.
For any $f\!\in\!\hbox{Hom}(\De^k,X)$, let
\begin{equation}\label{homot_e5}
D_Xf=f_{\#}D_{\De^k}\Id_{\De^k}.
\end{equation}
This construction defines a natural transformation $S_k\!\lra\!S_{k+1}$.
Since $D_{\De^k}(\Id_{\De^k})$ is a linear combination of linear maps, 
it is clear that the requirement (i) above is satisfied;
it remains to check (ii) and~(iii).\\

\noindent
(3) Suppose $f\!\in\!\hbox{Hom}(\De^k,X)$ and $\tau\!\in\!{\cal S}_k$, let
$s\!=\!f_{\#}\ti\tau\!\in\!S_k'(X)$.
By \e_ref{homot_e5}, \e_ref{homot_e4}, \e_ref{homot_e2}, and 
naturality~of~$D_X|_{S_{k-1}}$,
\begin{equation}\label{homot_e6}\begin{split}
D_X(f\!\circ\!\tau)
&=f_{\#}\tau_{\#}D_{\De^k}\Id_{\De^k}
=f_{\#}\tau_{\#}P_k
\big(\Id_{\De^k}+(-1)^{k+1}D_{\De^k}\partial_{\De^k}\Id_{\De^k}\big)\\
&=f_{\#}P_k\big(\tau+(-1)^{k+1}\tau_{\#}
                 D_{\De^k}\partial_{\De^k} \Id_{\De^k}\big)\\
&=f_{\#}P_k\big(\tau+(-1)^{k+1}D_{\De^k}\partial_{\De^k}\tau\big).
\end{split}\end{equation}
Thus,
\begin{equation}\label{homot_e7}
D_Xs=f_{\#}P_k\big(\ti\tau
+(-1)^{k+1}D_{\De^k}\partial_{\De^k}\ti\tau\big).
\end{equation}
By Lemma~\ref{complex}, the induction assumption~(ii), and~\e_ref{homot_e1b}, 
$S_k'(\De^k)$ is mapped into $S_*'(\De^k)$ 
by $D_{\De^k}\partial_{\De^k}$ and by~$P_k$.
Thus, by~\e_ref{homot_e7}, $D_X$ maps $S_k'(X)$ into~$S_{k+1}'(X)$.
Finally, by~\e_ref{homot_e7}, \e_ref{homot_e3}, and the inductive assumption~(iii),
\begin{equation*}\begin{split}
\partial_X D_Xs
&=\partial_X f_{\#}P_k\big(\ti\tau
\!+\!(-1)^{k+1}D_{\De^k}\partial_{\De^k}\ti\tau\big)
=f_{\#}\partial_{\De^k}P_k\ti\tau
+(-1)^{k+1}f_{\#}\partial_{\De^k}P_kD_{\De^k}\partial_{\De^k}\ti\tau\\
&=f_{\#}\big((-1)^{k+1}\ti\tau+P_k\partial_{\De^k}\tilde{\tau}\big)
+(-1)^{k+1}f_{\#}\big((-1)^{k+1}D_{\De^k}\partial_{\De^k}\tilde{\tau}
+P_k\partial_{\De^k}D_{\De^k}\partial_{\De^k}\tilde{\tau}\big)\\
&=\big((-1)^{k+1}s\!+\!D_X\partial_Xs\big)
+f_{\#}P_k\partial_{\De^k}\tilde{\tau}
+(-1)^{k+1}f_{\#}P_k\big((-1)^k\partial_{\De^k}\tilde{\tau}
\!+\!D_{\De^k}\partial_{\De^k}^2\tilde{\tau}\big)\\
&=(-1)^{k+1}s+D_X\partial_Xs.
\end{split}\end{equation*}

\begin{crl}
\label{acyclic_crl}
All homology groups of the complex $\big(S_*'(X),\partial_X|_{S_*'(X)}\big)$
are zero.
\end{crl}

\noindent
Let $\bar{S}_*(X)\!=\!S_*(X)/S_*'(X)$ and denote by 
$$\pi\!: S_*(X)\lra\bar{S}_*(X)$$ 
the projection map. 
Let $\bar{\partial}_X$ be boundary map on $\bar{S}_*(X)$
induced by~$\partial_X$.
We denote by $\bar{H}_*(X;\Bbb{Z})$ the homology groups of 
$\big(\bar{S}_*(X),\bar{\partial}_X\big)$.

\begin{prp}
\label{isom1_prp}
If $X$ is a topological space, the projection map 
$\pi\!:S_*(X)\!\lra\!\bar{S}_*(X)$ 
induces a natural isomorphism 
$H_*(X;\Bbb{Z})\!\lra\!\bar{H}_*(X;\Bbb{Z})$.
This isomorphism can be extended to relative homologies
to give an isomorphism of homology theories.
\end{prp}

\noindent
The first statement follows from the long exact sequence in homology
for the short exact sequence of chain complexes
$$0\lra S_*'(X)\lra S_*(X)\stackrel{\pi}{\lra}\bar{S}_*(X)\lra 0$$
and Corollary~\ref{acyclic_crl}.
The second statement follows from the first and the Five Lemma; 
see Lemma~24.3 in~\cite{Mu2}.\\

\noindent
For a simplicial complex $K$, there is a natural chain map from 
the {\it ordered} simplicial complex $C'_*(K)$ to the singular chain complex $S_*(|K|)$,
which induces isomorphism in homology.
If the vertices of $K$ are ordered, there is also a chain map
from $C'_*(K)$ to the {\it oriented} chain complex $C_*(K)$,
which induces a natural isomorphism in homology.
However, the chain map itself depends on the ordering of the vertices; 
see Section~34 in~\cite{Mu2}.
The advantage of the complex $\bar{S}_*(K)$ is that
there is a natural chain map from $C_*(K)$ to $\bar{S}_*(K)$,
which induces isomorphism in homology; this chain map is induced
by the natural chain map from $C'_*(K)$ to $S_*(|K|)$
described in Section~34 of~\cite{Mu2}. \\

\noindent
If $(X,\Bd X)$ is a compact oriented $n$-manifold, $(K,K',\eta)$ 
a triangulation of $(X,\Bd X)$, and for each 
$n$-dimensional simplex $\si\!\in\!K$, 
$$l_{\si}\!: \De^n\lra\si$$
is a linear
map such that $\eta\!\circ\!l_{\si}$ is orientation-preserving,
then the fundamental homology class $[X]\!\in\!H_n(X,\Bd X)$
is represented in $\bar{S}_k(X,\Bd X)$ by 
$$\sum_{\si\in K,\dim\si=n}\!\!\!\!\!\!\!\{\eta\!\circ\!l_{\si}\}\equiv
\sum_{\si\in K,\dim\si=n}\!\!\!\!\!\!\!\!\!\pi(\eta\!\circ\!l_{\si}),$$
where $\pi$ is as before. 
Note that 
$$\sum_{\si\in K,\dim\si=n}\!\!\!\!\!\!\!\!\eta\!\circ\!l_{\si}$$
may not even be a cycle in $S_k(X,\Bd X)$. 
It is definitely {\it not} a cycle if $\Bd X\!=\!\eset$ and $n$ is an even positive integer, 
as the boundary of each term $\eta\!\circ\!l_{\si}$ contains one more term with
coefficient~$+1$ than~$-1$.
Similarly, if 
$$h\!: (X,\Bd X)\lra(M,U)$$ 
is a continuous map, $h_*([X])\!\in\!H_k(M,U)$ is represented in $\bar{S}_k(M,U)$
by 
$$\sum_{\si\in K,\dim\si=n}\!\!\!\!\!\!\!\{h\!\circ\!\eta\!\circ\!l_{\si}\}.$$
Once again, the obvious preimage under $\pi$ of the above chain in
$S_k(M,U)$ may not be even a cycle.

\subsection{Combinatorics of Oriented Singular Homology}
\label{comb_subs}

\noindent
In this subsection we characterize cycles and boundaries in $\bar{S}_*(X)$
in a manner suitable for converting them to pseudocycles and pseudocycle 
equivalences in Subsection~\ref{psi_sec}.
We will use the two lemmas below to glue maps from standard simplices
together to construct smooth maps from smooth manifolds.\\

\noindent
The homology groups of a smooth manifolds~$X$ can be defined 
with the space $\hbox{Hom}(\De^k,X)$ of continuous maps from~$\De^k$ to~$X$
replaced by the space $C^{\i}(\De^k,X)$ of smooth maps;
this is a standard fact in differential topology. 
Note that the operator $D_X$ of Lemma~\ref{homot} maps smooth maps 
into linear combinations of smooth maps. 
Thus, all of the constructions of Subsection~\ref{singular} 
go through for the chain complexes based on elements in $C^{\i}(\De^k,X)$
instead of $\hbox{Hom}(\De^k,X)$.
Below $\bar{S}_*(X)$ will refer to the quotient complex based on  
such maps.\\

\noindent
If $s\!=\!\sum\limits_{j=1}\limits^{j=N}f_j$, where  
$f_j\!:\De^k\!\lra\!X$ is a continuous map for each~$j$, let
$${\cal C}_s=\big\{(j,p)\!: j\!=\!1,\ldots,N;~ p\!=\!0,\ldots,k\big\}.$$

\begin{lmm}
\label{gluing_l1}
If $k\!\ge\!1$ and $s\equiv\sum\limits_{j=1}\limits^{j=N}f_j$ 
determines a cycle in~$\bar{S}_k(X)$,
there exist a subset ${\cal D}_s\!\subset\!{\cal C}_s\!\times\!{\cal C}_s$
disjoint from the diagonal and a~map 
$$\tau\!: \D_s\lra {\cal S}_{k-1}, \qquad 
\big((j_1,p_1),(j_2,p_2)\big)\lra\tau_{(j_1,p_1),(j_2,p_2)},$$ 
such that\\
${}\quad~$ (i) if $\big((j_1,p_1),(j_2,p_2)\big)\!\in\!{\cal D}_s$, then
$\big((j_2,p_2),(j_1,p_1)\big)\!\in\!{\cal D}_s$;\\
${}~~~$ (ii) the projection ${\cal D}_s\!\lra\!{\cal C}_s$ on either coordinate
is a bijection;\\
${}~~\,$ (iii) for all $\big((j_1,p_1),(j_2,p_2)\big)\!\in\!{\cal D}_s$, 
\begin{gather}\label{psi_e1a}
\tau_{(j_2,p_2),(j_1,p_1)}=\tau_{(j_1,p_1),(j_2,p_2)}^{~-1}, \qquad
f_{j_2}\circ\io_{k,p_2}=
f_{j_1}\circ\io_{k,p_1}\circ\tau_{(j_1,p_1),(j_2,p_2)},\\
\label{psi_e1b}
\hbox{and}\qquad \sign\tau_{(j_1,p_1),(j_2,p_2)}=-(-1)^{p_1+p_2}.
\end{gather}\\
\end{lmm}

\noindent
This lemma follows from the assumption that $\bar\partial\{s\}\!=\!0$
and from the definition of $\bar{S}_*(X)$ in Subsection~\ref{singular}.
The terms appearing in the boundary of~$s$ are indexed by the set $\C_s$,
and the coefficient of the $(j,p)$th term is~$(-1)^p$.
Since~$s$ determines a cycle in~$\bar{S}_*(X)$,
these terms cancel in pairs, possibly after composition with an element 
$\tau\!\in\!{\cal S}_{k-1}$ and multiplying by~$\sign\tau$.
This operation does not change the equivalence class of a 
$(k\!-\!1)$-simplex in~$\bar{S}_{k-1}(X)$.

\begin{lmm}
\label{gluing_l2}
Suppose $k\!\ge\!1$, 
$$s_0\equiv\sum\limits_{j=1}\limits^{j=N_0}\{f_{0,j}\}, \qquad
s_1\equiv\sum\limits_{j=1}\limits^{j=N_1}\{f_{1,j}\}, \qquad
\ti{s}\equiv\sum\limits_{j=1}\limits^{j=\ti{N}}\ti{f}_j, \quad\hbox{and}\quad
\bar\partial\{\ti{s}\}=\{s_1\}-\{s_0\}\in\bar{S}_k(X).$$
Then there exist a subset $\D_{\ti{s}}\!\subset\!\C_{\ti{s}}\!\times\!\C_{\ti{s}}$
disjoint from the diagonal, subsets
${\cal C}_{\ti{s}}^{(0)},{\cal C}_{\ti{s}}^{(1)}\!\subset\!{\cal C}_{\ti{s}}$, 
and maps
\begin{gather*}
\ti\tau\!: \D_{\ti{s}}\lra {\cal S}_k, \quad 
\big((j_1,p_1),(j_2,p_2)\big)\lra\ti\tau_{(j_1,p_1),(j_2,p_2)},\\
(\ti{j}_i,\ti{p}_i)\!:\big\{1,\ldots,N_i\big\}\lra{\cal C}_{\ti{s}}^{(i)},
~~~\hbox{and}~~~
\ti\tau_i\!:\big\{1,\ldots,N_i\big\}\lra{\cal S}_k, ~~j\lra\ti\tau_{(i,j)}, \quad i=0,1,
\end{gather*}
such that\\
${}\quad~$ (i) if $\big((j_1,p_1),(j_2,p_2)\big)\!\in\!{\cal D}_{\tilde{s}}$, then
$\big((j_2,p_2),(j_1,p_1)\big)\!\in\!{\cal D}_{\tilde{s}}$;\\
${}~~~~$ (ii) the projection $\D_{\ti{s}}\!\lra\!\C_{\ti{s}}$  on either coordinate
is a bijection onto the complement of~$\C_{\ti{s}}^{(0)}\cup\C_{\ti{s}}^{(1)}$;\\
${}~~$ (iii) for all $\big((j_1,p_1),(j_2,p_2)\big)\!\in\!\D_{\ti{s}}$, 
\begin{gather}\label{psi_e2a}
\ti\tau_{(j_2,p_2),(j_1,p_1)}=\ti\tau_{(j_1,p_1),(j_2,p_2)}^{~-1}, \qquad
\ti{f}_{j_2}\circ\io_{k+1,p_2}=
\ti{f}_{j_1}\circ\io_{k+1,p_1}\circ\ti\tau_{(j_1,p_1),(j_2,p_2)},\\
\label{psi_e2b}
\hbox{and}\qquad \sign\ti\tau_{(j_1,p_1),(j_2,p_2)}=-(-1)^{p_1+p_2};
\end{gather}
${}~~~$  (iv) for all $i\!=\!0,1$ and $j\!=\!1,\ldots,N_i$, 
\begin{equation}\label{psi_e3}
\ti{f}_{\ti{j}_i(j)}\circ\io_{k+1,\ti{p}_i(j)}\circ\ti\tau_{(i,j)}=f_{i,j}
\qquad\hbox{and}\qquad \sign\ti\tau_{(i,j)}=-(-1)^{i+\ti{p}_i(j)};
\end{equation}
${}~~~~$  (v) $(\ti{j}_i,\ti{p}_i)$ is a bijection onto 
$\C_{\ti{s}}^{(i)}$ for $i\!=\!0,1$.
\end{lmm}

\noindent
This lemma follows from the assumption that 
$$\bar\partial\{\ti{s}\}=\{s_1\}-\{s_0\}.$$ 
The terms making up $\partial\ti{s}$ are indexed by the set~$\C_{\ti{s}}$.
By definition of~$\bar{S}_*(X)$, there exist disjoint subsets  
$\C_{\ti{s}}^{(0)}$ and $\C_{\ti{s}}^{(1)}$ of~$\C_{\ti{s}}$
such that for each $(j,p)\!\in\!\C_{\ti{s}}^{(1)}$
the $(j,p)$th term of $\partial\ti{s}$ equals one of the terms of~$s_i$,
after a composition with some $\ti\tau\!\in\!{\cal S}_k$ and multiplying
by $-(-1)^i\sign\ti\tau$.
The remaining terms of  $\C_{\ti{s}}$ must cancel in pairs, 
as in the case of Lemma~\ref{gluing_l1}.

\section{Integral Homology and  Pseudocycles}

\subsection{From Integral Cycles to Pseudocycles}
\label{psi_sec}

\noindent
In this subsection, we prove

\begin{prp}
\label{psi_prp}
If $X$ is a smooth manifold, there exists a homomorphism
$$\Psi_*\!: H_*(X;\Bbb{Z})\!\lra\!{\cal H}_*(X),$$
which is natural with respect to smooth~maps.
\end{prp}

\noindent
In the proof of Lemma~\ref{constr1}, we construct a homomorphism from 
the subgroup of cycles in~$\bar{S}_*(X)$ to~${\cal H}_*(X)$.
Starting with a cycle~$\{s\}$ as in Lemma~\ref{gluing_l1}, 
we will glue the functions~$f_j\!\circ\!\vp_k$ together,
where $\vp_k$ is the self-map of $\De^k$ provided by Lemma~\ref{bdpush_lmm}. 
These functions continue to satisfy the second equation in~\e_ref{psi_e1a},
i.e.
\begin{equation}\label{pertmaps_e1}
f_{j_2}\!\circ\!\vp_k\circ\io_{k,p_2}
=f_{j_1}\!\circ\!\vp_k\circ\io_{k,p_1}\circ\tau_{(j_1,p_1),(j_2,p_2)}
\qquad\forall~\big((j_1,p_1),(j_2,p_2)\big)\!\in\!\D_s,
\end{equation}
because $\vp_k\!=\!\id$ on $\De^k\!-\!\Int\,\De^k$ by the first equation 
in~\e_ref{bdpush_e1}.
The proof of Lemma~\ref{constr1} implements a construction suggested 
in Section~7.1 of~\cite{McSa}.\\

\noindent
Lemma~\ref{equiv1} shows that the map of Lemma~\ref{constr1}
descends to the homology groups. 
Starting with a chain~$\{\ti{s}\}$ as in Lemma~\ref{gluing_l2}, 
we will glue the functions $\ti{f}_j\!\circ\!\ti\vp_{k+1}\!\circ\!\vp_{k+1}$ together,
where $\ti\vp_{k+1}$ and $\vp_{k+1}$ are the self-maps of 
$\De^{k+1}$ provided by Lemma~\ref{bdpush_lmm}.
If $i\!=\!0,1$ and $j\!=\!1,\ldots,N_i$,
by the third equation in~\e_ref{bdpush_e2}, 
the second equation in~\e_ref{bdpush_e1},
and the first equation in~\e_ref{psi_e3}
\begin{equation*}\begin{split}
\ti{f}_{\ti{j}_i(j)}\!\circ\!\ti\vp_{k+1}\circ
\io_{k+1,\ti{p}_i(j)}\!\circ\!\ti\tau_{(i,j)}
=\ti{f}_{\ti{j}_i(j)}\circ
\io_{k+1,\ti{p}_i(j)}\!\circ\!\vp_k\!\circ\!\ti\tau_{(i,j)}
&=\ti{f}_{\ti{j}_i(j)}\circ
\io_{k+1,\ti{p}_i(j)}\!\circ\!\ti\tau_{(i,j)}\!\circ\!\vp_k\\
&=f_{i,j}\!\circ\!\vp_k.
\end{split}\end{equation*}
Since $\vp_{k+1}\!=\!\id$ on $\De^{k+1}\!-\!\Int\,\De^{k+1}$, it follows that
\begin{equation}\label{pertmaps_e2}
\ti{f}_{\ti{j}_i(j)}\!\circ\!\ti\vp_{k+1}\!\circ\!\vp_{k+1}\circ
\io_{k+1,\ti{p}_i(j)}\!\circ\!\ti\tau_{(i,j)}=f_{i,j}\!\circ\!\vp_k
\qquad\forall~j\!=\!1,\ldots,N_i,~i\!=\!0,1.
\end{equation}
Similarly, if $((j_1,p_1),(j_2,p_2))\!\in\!\D_{\ti{s}}$,
by the third equation in~\e_ref{bdpush_e2} used twice, 
the second equation in~\e_ref{psi_e2a},
and the second equation in~\e_ref{bdpush_e1},
\begin{equation*}\begin{split}
\ti{f}_{j_2}\!\circ\!\ti\vp_{k+1}\circ\io_{k+1,p_2}
&=\ti{f}_{j_2}\!\circ\!\io_{k+1,p_2}\!\circ\!\vp_k
=\ti{f}_{j_1}\!\circ\!\io_{k+1,p_1}\!\circ\!\ti\tau_{(j_1,p_1),(j_2,p_2)}
\circ\vp_k\\
&=\ti{f}_{j_1}\circ\io_{k+1,p_1}\!\circ\!\vp_k\circ\ti\tau_{(j_1,p_1),(j_2,p_2)}
=\ti{f}_{j_1}\!\circ\!\ti\vp_{k+1}\circ\io_{k+1,p_1}\!\circ\!\ti\tau_{(j_1,p_1),(j_2,p_2)}.
\end{split}\end{equation*}
Since $\vp_{k+1}\!=\!\id$ on $\De^{k+1}\!-\!\Int\,\De^{k+1}$, it follows that
\begin{equation}\label{pertmaps_e3}
\ti{f}_{j_2}\!\circ\!\ti\vp_{k+1}\!\circ\!\vp_{k+1}\circ\io_{k+1,p_2}
=\ti{f}_{j_1}\!\circ\!\ti\vp_{k+1}\!\circ\!\vp_{k+1}
\circ\io_{k+1,p_1}\!\circ\!\ti\tau_{(j_1,p_1),(j_2,p_2)}
\quad\forall~((j_1,p_1),(j_2,p_2))\!\in\!\D_{\ti{s}}.
\end{equation}
Thus, the functions $\ti{f}_j\!\circ\!\ti\vp_{k+1}\!\circ\!\vp_{k+1}$ are 
the analogues (in the sense of~Lemma~\ref{gluing_l2}) of the functions~$\ti{f}_j$
for the maps $f_{0,j}\!\circ\!\vp_k$ and~$f_{1,j}\!\circ\!\vp_k$.

\begin{lmm}
\label{constr1}
If $X$ is a smooth manifold,
every integral $k$-cycle in $X$, based on $C^{\i}(\De^k;X)$, 
determines an element of~${\cal H}_k(X)$.
\end{lmm}

\noindent
{\it Proof:} (1) If $k\!=\!0$, this is obvious. 
Suppose $k\!\ge\!1$ and 
$$s\equiv\sum_{j=1}^{j=N}f_j$$
determines a cycle in $\bar{S}_k(X)$.
Let ${\cal D}_s$ be the set provided by Lemma~\ref{gluing_l1}
and let $\tau\!:{\cal D}_s\!\lra\!{\cal S}_{k-1}$ be
the corresponding map.
Let 
\begin{gather*}
M'=\big(\bigsqcup_{j=1}^{j=N}\{j\}\!\times\!\De^k\big)\big/\sim,
\qquad\hbox{where}\\ 
\big(j_1,\io_{k,p_1}(\tau_{(j_1,p_1),(j_2,p_2)}(t))\big)
\!\sim\!\big(j_2,\io_{k,p_2}(t)\big)
\quad\forall~
\big((j_1,p_1),(j_2,p_2)\big)\!\in\!{\cal D}_s,~t\!\in\!\De^{k-1}.
\end{gather*}
Let $\pi$ be the quotient map.
Define 
\begin{equation}\label{fmapdfn_e}
F\!: M'\lra X \qquad\hbox{by}\quad F\big([j,t]\big)= f_j\big(\vp_k(t)\big).
\end{equation}
This map is well-defined by~\e_ref{pertmaps_e1} and continuous by 
the universal property of the quotient topology; see Theorem~22.2 in~\cite{Mu1}.
Let $M$ be the complement in $M'$ of the~set 
$$\pi\Big(\bigsqcup_{j=1}^{j=N}\{j\}\!\times\!Y\Big),$$
where $Y$ is the $(k\!-\!2)$-skeleton of $\De^k$.
By continuity of~$F$, compactness of $M'$, and the first equation in~\e_ref{bdpush_e1},
\begin{equation}\label{constr1_e3}
\Bd F|_M=F(M'\!-\!M)= \bigcup_{j=1}^{j=N}f_j\big(\vp_k(Y)\big)
= \bigcup_{j=1}^{j=N}f_j(Y).
\end{equation}
Since $f_j|_{\Int\si}$ is smooth for all $j\!=\!1,\ldots,N$ and 
all simplices $\si\!\subset\!\De^k$,
$\Bd F|_M$ has dimension at most $k\!-\!2$ by~\e_ref{constr1_e3}.
Thus, $F|_M$ is a $k$-pseudocycle, provided 
$M$ is a smooth oriented manifold and $F|_M$ is a smooth~map.
This is shown in~(2) below.\\

\noindent
(2) Let $[j,t]\!\in\!M$ be any point. 
If $t\!\in\!\Int\De^k$, then 
 $\pi(\{j\}\!\times\!\Int\De^k)$ is an open set about~$[j,t]$, 
which is naturally homeomorphic to~$\Int\De^k$. 
If 
$$[j,t]=\big[j_1,\io_{k,p_1}(t_1)\big]=\big[j_2,\io_{k,p_2}(t_2)\big]$$
with $(j_1,p_1)\!\neq\!(j_2,p_2)$ and $t_1\!\in\!\Int\,\De^{k-1}$, let 
$$U=\pi\big(\{j_1\}\!\times\!U^k_{p_1}\big)
\cup\pi\big(\{j_2\}\!\times\!U^k_{p_2}\big).$$
This is an open neighborhood of $[j,t]$ in $M$.
It is homeomorphic in a canonical way to the disjoint union~of 
$U^k_{p_1}$ and $U^k_{p_2}$ with $\Int\De^k_{p_1}\!\subset\!U^k_{p_1}$ and 
$\Int\De^k_{p_2}\!\subset\!U^k_{p_2}$ identified by the linear map
\begin{equation}\label{constr1_e3a}
\io_{k,p_1}\circ\tau_{(j_1,p_1),(j_2,p_2)}\circ\io_{k,p_2}^{-1}\!:
\Int\De^k_{p_2}\lra\Int\De^k_{p_1}
\end{equation}
and thus to an open subset of $\R^k$.
By~\e_ref{psi_e1b}, the transition map~\e_ref{constr1_e3a} 
is orientation-reversing if
the open simplices $\Int\De^k_{p_1}$ and  $\Int\De^k_{p_2}$ are oriented 
as boundaries of the $k$-manifolds $U^k_{p_1}$ and $U^k_{p_2}$
with their natural orientations. 
This means that the induced orientations of~$T_pU$ 
coming from the two $k$-manifolds with boundary agree.
On any nonempty overlap of this coordinate chart with any other
coordinate chart, the transition map is the identity map on an open subset
of $\Int\De^k$. 
Thus, $M$ is a smooth oriented manifold. 
The map~$F$ is smooth on $\{j\}\!\times\!\Int\De^k$ 
for all~$j$ by our assumptions on~$F$.
If 
$$[j,t]=\big[j_1,\io_{k,p_1}(t_1)\big]=\big[j_2,\io_{k,p_2}(t_2)\big],$$
then $F$ is smooth on the open set $U$, defined as above, because
it is smooth on 
$$\pi\big(\{j_1\}\!\times\!U^k_{p_1}\big) \quad\hbox{and}\quad
\pi\big(\{j_2\}\!\times\!U^k_{p_2}\big),$$  
and all derivatives in the direction normal to 
$\pi(\{j_1\}\!\times\!\Int\De^k_{p_1})$ vanish 
by the first equation in~\e_ref{bdpush_e1}.\\

\noindent
{\it Remark:} The pseudocycle $F|_M$ constructed above depends on the choice
of~$\D_s$ and~$\tau$.
However, as the next lemma shows, the image of~$F|_M$
in ${\cal H}_k(X)$ depends only on~$[\{s\}]$.

\begin{lmm}
\label{equiv1}
Under the construction of Lemma~\ref{constr1}, 
homologous $k$-cycles determine
the same equivalence class of pseudocycles in ${\cal H}_k(X)$.
\end{lmm}

\noindent
{\it Proof:} 
(1) If $k\!=\!0$, this is obvious. 
Suppose $k\!>\!0$ and 
$$s_0\equiv\sum_{j=1}^{j=N_0}f_{0,j} \qquad\hbox{and}\qquad 
s_1\equiv\sum_{j=1}^{j=N_1}f_{1,j}$$
determine two homologous $k$-cycles in $\bar{S}_k(X)$.
Let ${\cal D}_{s_0}$ and ${\cal D}_{s_1}$ be the sets provided
by Lemma~\ref{gluing_l1} and let $\tau_0$ and
$\tau_1$ be the corresponding maps into~${\cal S}_{k-1}$.
Denote by $(M_0',M_0,F_0)$ and $(M_1',M_1,F_1)$
the triples constructed in the proof of Lemma~\ref{constr1}
corresponding to $s_0$ and~$s_1$.
Let 
$$\ti{s}=\sum_{j=1}^{j=\ti{N}}\ti{f}_j\in S_{k+1}(X)$$
be such that 
$$\bar\partial\{\ti{s}\}=\{s_1\}-\{s_0\}\in\bar{S}_k(X).$$
Denote by  ${\cal C}_{\ti{s}}^{(0)}$, ${\cal C}_{\ti{s}}^{(1)}$, 
${\cal D}_{\ti{s}}$, $(\ti{j}_i,\ti{p}_i,\ti\tau_i)$, and $\ti\tau$
the corresponding objects of Lemma~\ref{gluing_l2}.\\

\noindent
{\it Remark:}
Proceeding as in the proof of Lemma~\ref{constr1},
we can turn $\ti{s}$ into a pseudocycle equivalence $(\ti{M}^*,\ti{F})$ 
between two pseudocycles $(M_0^*,F_0)$ and $(M^*,F_1)$
by gluing across codimension-one faces.
Unfortunately, $M_0^*$ and $M_1^*$ are not the entire manifolds $M_0$ and $M_1$;
they are missing the $(k\!-\!1)$-simplices of $M_0$ and~$M_1$.
This issue is resolved in~(2) below by adding collars to $\ti{M}^*$: 
$(n\!+\!1)$-manifolds that begin with $M_i^*$ and end with~$M_i^*$.\\

\noindent
(2) Let $I\!=\![0,1]$. Put 
\begin{gather*}
\tilde{M}'=\Big(\bigsqcup_{j=1}^{j=\tilde{N}}\{j\}\!\times\!\De^{k+1} \sqcup
\bigsqcup_{i=0,1}\!\{i\}\!\times\!I\!\times\!M_i'\Big)
\Big/\sim,
\qquad\hbox{where}\\ 
\big(j_1,\io_{k+1,p_1}(\ti\tau_{(j_1,p_1),(j_2,p_2)}(t))\big)
\sim\big(j_2,\io_{k+1,p_2}(t)\big)
\quad\forall~
\big((j_1,p_1),(j_2,p_2)\big)\!\in\!\ti\D_{\ti{s}},~t\!\in\!\De^k,\\
\big(i,1\!-\!i,\pi(j,t)\big)\sim
\big(\ti{j}_i(j),\io_{k+1,\ti{p}_i(j)}(\ti\tau_{i,j}(t))\big)
\quad\forall~t\!\in\!\De^k,~j\!=\!1,\ldots,N_i,~i\!=\!0,1.
\end{gather*}
Let 
$$\ti\pi\!: \bigsqcup_{j=1}^{j=\tilde{N}}\{j\}\!\times\!\De^{k+1}
\sqcup \bigsqcup_{i=0,1}\!\{i\}\!\times\!I\!\times\!M_i' \lra \ti{M}'$$
be the quotient map.
Define 
\begin{equation*}
\ti{F}\!:\ti{M}'\lra X
\qquad\hbox{by}\qquad
\begin{aligned}
&\ti{F}\big([j,t]\big)=\ti{f}_j\big(\ti\vp_{k+1}(\vp_{k+1}(t))\big) 
\quad&\forall~t\!\in\!\De^{k+1},~j\!=\!1,\ldots,\ti{N};\\
&\ti{F}\big([i,s,x]\big)=F_i(x)\quad&\forall~s\!\in\!I,~x\!\in\!M_i',~i\!=\!0,1.
\end{aligned}
\end{equation*}
This map is well-defined by \e_ref{pertmaps_e2}, 
\e_ref{pertmaps_e3}, and~\e_ref{fmapdfn_e} and is continuous by the universal property
of the quotient topology.
Let $\ti{M}$ be the complement in~$\ti{M}'$ of the~set
$$\ti\pi\Big(\bigsqcup_{j=1}^{j=\ti{N}}\{j\}\!\times\!\ti{Y}
\sqcup\bigsqcup_{i=0,1}\!\{i\}\!\times\!I\!\times\!(\ti{M}_i'\!-\!\ti{M})\Big),$$
where $\ti{Y}$ is the $(k\!-\!1)$-skeleton of $\De^{k+1}$.
By continuity of~$\ti{F}$, compactness of $\ti{M}'$, and the first equation in~\e_ref{bdpush_e2},
\begin{equation}\label{equiv1_e3}
\Bd \ti{F}|_{\ti{M}}=\ti{F}(\ti{M}'\!-\!\ti{M})
= \bigcup_{j=1}^{j=\ti{N}}
\ti{f}_j\big(\ti\vp_{k+1}(\vp_{k+1}(\ti{Y}))\big)
\cup\bigcup_{i=0,1}\!f_{i,j}\big(\vp_k(Y)\big)\
=\bigcup_{j=1}^{j=\ti{N}}\ti{f}_j\big(\ti{Y}\big).
\end{equation}
Since $\ti{f}_j|_{\Int\si}$ is smooth for all $j\!=\!1,\ldots,\ti{N}$ and 
all simplices $\si\!\subset\!\De^{k+1}$,
$\Bd \ti{F}|_{\ti{M}}$ has dimension at most $k\!-\!1$ by~\e_ref{equiv1_e3}.
Thus, $\ti{F}|_{\ti{M}}$ is a pseudocycle equivalence between $F_0|_{M_0}$
and $F_1|_{M_1}$, provided $\ti{M}$ is a smooth oriented manifold, 
$\ti{F}|_{\ti{M}}$ is a smooth~map, and
$$\partial \big(\ti{F}|_{\ti{M}}\big)=F_1|_{M_1}-F_0|_{M_0}.$$
This is shown in~(3) below.\\

\noindent
(3) By the proof of Lemma~\ref{constr1}, for each $i\!=\!0,1$,
$$\ti{M}_i\equiv\ti{M}\cap \ti\pi\big(\{i\}\!\times\!I\!\times\!M_i'\big)
\approx I\!\times\!M_i-\{1\!-\!i\}\!\times\!\bigcup_{j=1}^{j=N_i}\!
\pi\big(\{j\}\!\times\!(\De^k\!-\!\Int\De^k)\big)$$
is a smooth oriented manifold with boundary,
$$\partial\ti{M}_i
\approx(-1)^{i+1}\big(\{i\}\!\times\!M_i\big)
+(-1)^i\Big(\{1\!-\!i\}\!\times\!\!\bigcup_{j=1}^{j=N_i}\{j\}\!\times\!\Int\De^k\Big),$$
and $\ti{F}$ restricts to a smooth map on $\ti{M}_i$.
By the same argument as in~(2) of the proof of Lemma~\ref{constr1},
the space 
$$\ti{M}^*\equiv
\ti\pi\big(\bigsqcup_{j=1}^{j=\ti{N}}\!\{j\}\!\times\!(\De^{k+1}\!-\!\ti{Y})\big)$$
is a smooth oriented manifold with boundary,
$$\partial\ti{M}^*
\approx\bigsqcup_{i=0,1}\bigsqcup_{(j,p)\in\C_{\ti{s}}^{(i)}}\!\!\!\!\!
\{j\}\!\times\Int\De^{k+1}_p,$$
and $\ti{F}$ restricts to a smooth map on $\ti{M}^*$.\\

\noindent
The topological space $\ti{M}$ is obtained by identifying
the $(k\!+\!1)$-manifolds $\ti{M}^*$ and $\ti{M}_i$ along
components of their boundaries via the~map
$$\{1\!-\!i\}\!\times\bigcup_{j=1}^{j=N_i}\!\{j\}\!\times\!\Int\De^k
\lra\bigsqcup_{(j,p)\in\C_{\ti{s}}^{(i)}}\!\!\!\!\!
\{j\}\!\times\Int\De^{k+1}_p,\qquad
\big(1\!-\!i,j,t)\lra \big(\ti{j}_i(j),\io_{k+1,\ti{p}_i(j)}(t)\big).$$
By \e_ref{signch_e} and the second equation in~\e_ref{psi_e3}, 
this overlap map is orientation-reversing.
Therefore, similarly to~(2) in the proof of Lemma~\ref{constr1}, 
it follows that $\ti{M}$ is a smooth oriented manifold with boundary and
$$\partial\ti{M}=\ti\pi(\{1\}\!\times\!\{0\}\!\times\!M_1\big)
-\ti\pi(\{0\}\!\times\!\{1\}\!\times\!M_0\big)
\approx M_1-M_0.$$
It is immediate that under this identification,
$$\ti{F}|_{M_i}=F_i \qquad\forall~ i=0,1.$$
Finally, as in the proof of Lemma~\ref{constr1}, all derivatives of
the smooth maps  $\ti{F}|_{M^*}$ and $\ti{F}|_{M_i}$ in the normal directions 
to their boundaries vanish.
Therefore, $\ti{F}|_{\ti{M}}$ is smooth and is a pseudocycle equivalence 
from $F_0|_{M_0}$ to~$F_1|_{M_1}$, as claimed.

\subsection{From Pseudocycles to Integral Cycles}
\label{phi_sec}

\noindent
In this subsection, we prove

\begin{prp}
\label{phi_prp}
If $X$ is a smooth manifold, there exists a homomorphism 
$$\Phi_*\!: {\cal H}_*(X)\!\lra\!H_*(X;\Z),$$
which is natural with respect to smooth maps.
\end{prp}

\begin{lmm}
\label{euler}
Every $k$-pseudocycle determines a class in $H_k(X;\Bbb{Z})$.
\end{lmm}

\noindent
{\it Proof:} (1) Suppose $h\!: M\!\lra\!X$ is a $k$-pseudocycle and 
$f\!: N\!\lra\!X$ a smooth map such that 
$$\dim N=k\!-\!2 \qquad\hbox{and}\qquad \Bd h\subset\Im f.$$
By Proposition~\ref{neighb_prp}, there exists an open neighborhood $U$
of $\Bd h$ in~$X$ such that 
$$H_l(U;\Z)=0 \qquad\forall l>k\!-\!2.$$ 
Let $K\!=\!M\!-\!h^{-1}(U)$.
Since the closure of $h(M)$ is compact in $X$, 
$K$ is a compact subset of $M$ by definition of $\Bd h$. 
Let $V$ be an open neighborhood of $K$ in $M$ 
such that $\bar{V}$ is a compact manifold with boundary.
It inherits an orientation from the orientation of~$M$ and thus defines a homology 
$$[\bar{V}]\in H_k(\bar{V},\Bd\bar{V};\Z).$$
Put
\begin{equation}\label{euler_e3}
[h]=h_*([\bar{V}])\in H_k(X,U;\Z)\approx H_k(X;\Z),
\end{equation}
where 
\begin{equation}\label{euler_e4}
h_*\!: H_k(\bar{V},\Bd\bar{V};\Z)\!\lra\! H_k(X,U;\Z)
\end{equation}
is the homology homomorphism induced by $h$.
The isomorphism in~\e_ref{euler_e3} is induced by inclusion. 
It is an isomorphism by the assumption on the homology of~$U$ as follows
from the long exact sequence in homology for the pair~$(X,U)$.\\

\noindent 
(2) The homology class $[h]$ is independent of the choice of $V$.
Suppose $V'$ is another choice such that $\bar{V}\!\subset\!V'$.
Choose a triangulation of $\bar{V}'$ extending some triangulation
of $(\Bd\bar{V})\bigcup(\Bd\bar{V}')$; 
such a triangulation exists by Section~16 in~\cite{Mu2}. 
The cycles 
$$h_*([\bar{V}]),h_*([\bar{V}'])\in H_k(X,U;\Z)$$
then differ by singular simplices lying in $U$;
see discussion at the end of Subsection~\ref{singular}.
Thus, 
$$h_*([\bar{V}'])=h_*([\bar{V}])\in H_k(X,U;\Z).$$\\

\noindent
(3) The cycle $[h]$ is also independent of the choice of $U$. 
Suppose $U'\!\subset\!U$ is another choice. 
By~(2), it can be assumed that $V$ and $V'$ chosen as in~(1) are the same. 
Since the isomorphism in~\e_ref{euler_e3} is the composite of isomorphisms
$$H_k(X;\Z)\lra H_k(X,U';\Z)\lra H_k(X,U;\Z)$$ 
induced by inclusions and the homomorphism~\e_ref{euler_e4} is the composition
$$H_k(\bar{V},\Bd\bar{V};\Z)\lra H_k(X,U';\Z)\lra H_k(X,U;\Z),$$ 
the  homology classes obtained in $H_k(X;\Z)$ from $U$ and $U'$ are equal.
Finally, if $U$ and $U'$ are two arbitrary choices of open sets
in~(1), by Proposition~\ref{neighb_prp} there exists a third choice
$U''\!\subset\!U\!\cap\!U'$.

\begin{lmm}
Equivalent $k$-pseudocycles determine the same class in $H_k(X,\Z)$.
\end{lmm}

\noindent
{\it Proof:} Suppose $h_i\!: M_i\!\lra\!X$, $i\!=\!0,1$, are two equivalent
$k$-pseudocycles and $\ti{h}\!: \ti{M}\!\lra\!X$ is an equivalence between them. 
In particular, $\ti{M}$ is oriented, 
$$\partial\ti{M}=\!M_1\!-\!M_0, \qquad\hbox{and}\qquad \ti{h}|_{M_i}=h_i.$$ 
Let $\ti{U}$ be an open neighborhood of $\Bd\ti{h}$ in~X such that 
$$H_l(\ti{U};\Z)=0  \qquad\forall~l>k\!-\!1.$$
Let $U_i$ be  an open neighborhood of $\Bd h_i\!\subset\!\Bd\tilde{h}$ 
in $\tilde{U}$ such that 
$$H_l(U_i;\Bbb{Z})=0 \qquad\forall~l>k\!-\!2,$$ 
as provided by Proposition~\ref{neighb_prp}.
Let $V_i\!\subset\!M_i$ be a choice of an open set as in (1) 
of the proof of Lemma~\ref{euler}.
For $i\!=\!0,1$, choose a triangulation of $M_i$ that extends a triangulation
of $\Bd\bar{V}_i$.
Extend these two triangulations to a triangulation 
$\ti{T}\!=\!(\ti{K},\ti\eta)$ of~$\ti{M}$. 
Let $K$ be a finite sub-complex of~$\ti{K}$ such~that
$$V_0,V_1\subset\tilde{\eta}(|K|)  \qquad\hbox{and}\qquad
\ti{M}-\ti{h}^{-1}(\ti{U})\subset\ti\eta(\Int|K|).$$
Such a subcomplex exists because $\ti{h}(\ti{M})$ is a pre-compact subset of $X$
and thus $\ti{M}\!-\!\ti{h}^{-1}(\ti{U})$ is a compact subset of~$\ti{M}$.
Put
$$K_i=\big\{\si\!\in\!K\!: \eta(\si)\!\subset\!\bar{V}_i\big\} 
\qquad\hbox{for}~~i=0,1.$$
By the proof of Lemma~\ref{euler}, 
$(K_i,\ti{h}\!\circ\!\ti\eta|_{|K_i|})$ determines 
the homology class  $[h_i]\!\in\!H_k(X,U_i;\Z)$. 
Let $[h_i']$ denote its image in $H_k(X,\ti{U};\Z)$
under the homomorphism induced by inclusion.
The above assumptions on~$K$ imply that 
$$\partial(K,\ti{h}\circ\ti\eta|_K)=
(K_1,\ti{h}\circ\ti\eta|_{K_1})-(K_0,\ti{h}\circ\ti\eta|_{K_0})$$
in $\bar{S}(M,\ti{U})$. 
Thus,  
$$[h_0']=[h_1']\in H_k(X,\ti{U};\Z),$$
and this class lies in the image of the homomorphism 
\begin{equation}\label{equiv_e11}
H_k(X;\Z)\lra H_k(X,\ti{U};\Z)
\end{equation}
induced by inclusion. 
This map is equal to the composites
\begin{gather*}
H_k(X;\Z)\lra H_k(X,U_0;\Z)\lra H_k(X,\ti{U};\Z),\\
H_k(X;\Z)\lra H_k(X,U_1;\Z)\lra H_k(X,\ti{U};\Z).
\end{gather*}
Since $H_k(\ti{U};\Z)\!=\!0$, the homomorphism~\e_ref{equiv_e11}
is injective.
Thus, $[h_0]$ and $[h_1]$ come from the same element of $H_k(X;\Z)$.

\subsection{Isomorphism of Homology Theories}
\label{isom_sec}

\noindent
In this subsection we conclude the proof of Theorem~\ref{main_thm}.

\begin{lmm}
\label{isom_l1}
If $X$ is a smooth manifold, the composition
$$\Phi_*\circ \Psi_*\!: 
H_*(X;\Z)\lra{\cal H}_*(X) \lra H_*(X;\Z)$$
is the identity map on $H_*(X;\Z)$.
\end{lmm}

\noindent
{\it Proof:} Suppose
$$\{s\}=\sum_{j=1}^N\{f_j\}\in\bar{S}_k(X)$$ 
is a cycle and $F\!: M\!\lra\!X$ is a pseudocycle corresponding to~$s$
via the construction of Lemma~\ref{constr1}.
Recall that $M$ is the complement of the $(k\!-\!2)$-simplices in a compact space~$M'$
and $F$ is the restriction of a continuous map $F'\!: M'\!\lra\!X$ induced by the~maps
$$f_j\!\circ\vp_k\!:\De^k\lra X, \qquad j=1,\ldots,N.$$
Since $\vp_k$ is homotopic to the identity on $\De^k$, with boundary fixed,
\begin{equation}\label{isom_l1e3}
f_j\!\circ\vp_k-f_j\in\partial S_{k+1}(X) \qquad\forall~j=1,\ldots,N.
\end{equation}
Let $U$ be a neighborhood of $\Bd F$ such that
$$H_l(U;\Z)=0 \qquad\forall~ l>k\!-\!2.$$ 
Put $K\!=\!M\!-\!f^{-1}(\vp_k^{-1}(U))$. 
Let $V$ be a pre-compact neighborhood of $K$ such that $(\bar{V},\partial\bar{V})$
is a smooth manifold with boundary. 
Choose a triangulation $T\!=\!(K,\eta)$ of $(\bar{V},\partial\bar{V})$ 
such that every $k$-simplex of~$T$ is contained in a set of 
the form $\pi(\{j\}\!\times\!\De^k)$, where $\pi$ is as in the proof of Lemma~\ref{constr1}. 
Put
$$K_j=\big\{\si\!\in\!K\!: \eta(\si)\!\subset\!\pi(\{j\}\!\times\!\De^k)\big\},
\qquad
K_j^{\top}=\big\{\si\!\in\!K_j\!: \dim\si=k\big\}.$$
Let $\ti{T}_j\!=\!(\ti{K}_j,\it\eta_j)$ be a triangulation of a subset of
$\De^k$ that along with $K_j$ gives a triangulation of~$\De^k$. 
Put 
$$\ti{K}_j^{\top}=\big\{\si\!\in\!\ti{K}_j\!: \dim\si=k\big\}.$$
By definition of~$T$,
\begin{equation}\label{isom_l1e5}
f_j\!\circ\!\vp_k\big(\it\eta_j(\si)\big)\subset U
\qquad\forall~\si\in\ti{K}_j^{\top}.
\end{equation}
Furthermore, by~\e_ref{isom_l1e3}
\begin{equation}\label{isom_l1e7}\begin{split}
\{s\}&=\sum_{\si\in K^{\top}}\{f_j\!\circ\!\vp_k\circ\eta\circ l_{\si}\}\\
&=\sum_{j=1}^N\sum_{\si\in K_j^{\top}}\{f_j\!\circ\!\vp_k\circ\eta\circ l_{\si}\}
+\sum_{j=1}^N\sum_{\si\in\ti{K}_j^{\top}}\{f_j\!\circ\!\vp_k\circ\ti\eta_j\circ l_{\si}\}
\quad\mod \bar\partial \bar{S}_{k+1}(X), 
\end{split}\end{equation}
since subdivisions of cycles do not change the homology class.
By the proof of Lemma~\ref{euler}, the first sum on the right-hand side 
of~\e_ref{isom_l1e7} represents $[F]$ in $\bar{S}_k(X,U)$.
By~\e_ref{isom_l1e5},  the second sum lies in $\bar{S}_k(U)$.
Since the sum of the two terms is a cycle in $\bar{S}_k(X)$,
it must represent $[F]$ in $\bar{S}_k(X)$. Thus, 
$$\{F\}=\{s\}\in H_k(X;\Z),$$ 
and the claim~follows.

\begin{lmm}
\label{isom_l2}
If $X$ is a smooth manifold, the homomorphism
$\Phi_*\!: {\cal H}_*(X)\!\lra\!H_*(X;\Z)$ is injective.
\end{lmm}

\noindent
We will assume that a $k$-pseudocycle $h\!: M'\!\lra\!X$ determines
the zero homology class via the construction of Lemma~\ref{euler}
and show that it must be the boundary of a smooth map 
$\ti{F}\!:\ti{M}\!\lra\!X$ in the sense of pseudocycles.
A new difficulty here is that $M'$ need not be compact,
and therefore $\ti{M}$ may need to be constructed from infinitely many 
$(k\!+\!1)$-simplices.
This will be achieved as the limit of finite stages~$\ti{M}_i$,
so that as $i\!\in\!\Z^+$ increases $M'$ will be the pseudocycle boundary of
$\ti{M}_i$ ``modulo'' smaller and smaller neighborhoods~$U_i$ of~$\Bd h$.
As in the proof of Lemma~\ref{equiv1}, we will also need to attach a collar
to the $(k\!+\!1)$-manifold $\ti{M}^*$ obtained directly from a bounding chain.\\

\noindent
{\it Proof:}
(1) Suppose a $k$-pseudocycle $h\!: M'\!\lra\!X$ determines
the zero homology class. 
It can be assumed that $k\!\ge\!1$; otherwise, there is nothing to prove.
Let $\big\{U_i\big\}_{i=1}^{\i}$ be a sequence of open pre-compact neighborhoods of~$\Bd h$ in~$X$
such that 
$$U_{i+1}\subset U_i, \qquad \bigcap_{i=1}^{\i}U_i=\Bd h,
\quad\hbox{and}\quad
H_l(U_i;\Bbb{Z})=0 \quad\forall~ l>k\!-\!2.$$
Existence of such a collection follows from Proposition~\ref{neighb_prp}
and metrizability of any manifold.
Let $\big\{V_i\big\}_{i=1}^{\i}$ be a corresponding collection of open
sets in~$M'$ as in (1) of the proof of Lemma~\ref{euler}.
It can be assumed that $\bar{V}_i\!\subset\!V_{i+1}$.
Choose a triangulation $T\!=\!(K,\eta)$ of $M'$ that extends a triangulation
of $\bigcup\limits_{i=1}\limits^{\i}\Bd\bar{V}_i$. 
Let
$$K^{\top}=\big\{\si\!\in\!K\!:\dim\si\!=\!k\big\}, \qquad
\C_{\eta}=\big\{(\si,p)\!:\si\!\in\!K^{\top},~p\!=\!0,1,\ldots,k\big\}.$$
For each $\si\!\in\!K^{\top}$, let 
$$l_{\si}\!: \De^k\lra\si\subset |K|\subset\R^{\i}$$
be a linear map such that $\eta\!\circ\!l_{\si}$ is orientation-preserving. 
Put
\begin{gather*}
f_{\si}=h\circ\eta\circ l_{\si}\qquad\forall~\si\!\in\!K^{\top} \qquad\hbox{and}\\
\D_{\eta}=\big\{((\si_1,p_1),(\si_2,p_2))\!\in\!\C_{\eta}\!\times\!\C_{\eta}:
(\si_1,p_1)\!\neq\!(\si_2,p_2),~
l_{{\si}_1}(\De^k_{p_1})\!=\!l_{{\si}_2}(\De^k_{p_2})\big\}.
\end{gather*}
For each $((\si_1,p_1),(\si_2,p_2))\!\in\!\D_{\eta}$, define
$$\tau_{(\si_1,p_1),(\si_2,p_2)}\in{\cal S}_{k-1} \qquad\hbox{by}\qquad
l_{\si_2}\!\circ\!\io_{k,p_2}=
l_{\si_1}\!\circ\!\io_{k,p_1}\circ\tau_{(\si_1,p_1),(\si_2,p_2)}.$$
Since $K$ is an oriented simplicial complex,
$$\D_{\eta}\subset\C_{\eta}\!\times\!\C_{\eta} \qquad\hbox{and}\qquad
\tau\!:\D_{\eta}\lra{\cal S}_{k-1}$$
satisfy (i)-(iii) of Lemma~\ref{gluing_l1}.
Furthermore, $M'$ is the topological space corresponding to 
$(\C_{\eta},\D_{\eta},\tau)$ via the construction of Lemma~\ref{constr1}
and $h$ is the continuous map described by
$$h|_{\pi(\si\times\De^k)}=f_{\si}.$$
As in the proof of Lemma~\ref{constr1}, let $M$ be the complement of 
the $(k\!-\!2)$-simplices in $M'$; the pseudocycles $h$ and $h|_M$ are equivalent.
Since $\vp_k$ is homotopic to the identity on~$\De^k$ with boundary fixed,
the pseudocycle $h|_M$ is in turn equivalent to the pseudocycle $F|_M$,
where as in the proof of Lemma~\ref{constr1}
$$F\!:M'\lra X, \qquad F\circ \eta\circ l_{\si}=f_{\si}\!\circ\!\vp_k.$$\\

\noindent
(2) For each $i\!\ge\!1$, let
$$K_i^{\top}=\big\{\si\!\in\!K^{\top}\!: \eta(\si)\!\subset\!\bar{V}_i\big\},
\qquad \C_{\eta;i}=\big\{(\si,p)\!\in\!\C_{\eta}:\si\!\in\!K_i^{\top}\big\},
\quad\hbox{and}\quad
\D_{\eta;i}=\D_{\eta}\!\cap\!(\C_{\eta;i}\!\times\!\C_{\eta;i}).$$
By construction of $[h]$, for every $i\!\ge\!1$,
there exists a singular chain 
$$s_i\equiv\sum_{j=1}^{N_i}f_{i,j}\in S_k(U_i)$$ 
such that
$$\sum_{\si\in K_i^{\top}}\{h\circ\eta\circ l_{\si}\}+\{s_i\}$$
is a cycle in $\bar{S}_k(X)$ representing $[h]$. 
Similarly to Lemma~\ref{gluing_l1}, there exist a symmetric subset 
$${\cal D}_i\subset({\cal C}_{\eta;i}\!\sqcup\!{\cal C}_{s_i}) 
\!\times\! ({\cal C}_{\eta;i}\!\sqcup\!{\cal C}_{s_i})$$
disjoint from the diagonal and a map
$$\tau_i\!: {\cal D}_i\lra {\cal S}_{k-1}$$ 
such that\\
${}~~~~$ (i) ${\cal D}_{\eta;i}\!\subset\!{\cal D}_i$ and 
$\tau_i|_{{\cal D}_{\eta;i}}\!=\!\tau|_{{\cal D}_{\eta;i}}$;\\
${}~~~$ (ii) the projection map 
${\cal D}_i\lra{\cal C}_{\eta;i}\!\sqcup\!{\cal C}_{s_i}$
on either coordinate is a bijection;\\
${}~~$ (iii) for all $((j_1,p_1),(j_2,p_2))\!\in\!{\cal D}_i$, 
\begin{gather*}
\tau_{(j_2,p_2),(j_1,p_1)}=\tau_{(j_1,p_1),(j_2,p_2)}^{~-1}, \qquad
f_{i,j_2}\circ\io_{k,p_2}=
f_{i,j_1}\circ\io_{k,p_1}\circ\tau_{(j_1,p_1),(j_2,p_2)},\\
\hbox{and}\qquad \sign\tau_{(j_1,p_1),(j_2,p_2)}=-(-1)^{p_1+p_2},
\end{gather*}
where $f_{i,\si}\!\equiv\!f_{\si}$ for all $\si\!\in\!K_i^{\top}$.\\

\noindent
(3) By (2), for each $i\!\ge\!2$ 
$$\sum_{\si\in K_i^{\top}-K_{i-1}^{\top}}\!\!\!\!\!\!
\{h\circ\eta\circ l_{\si}\}+\{s_i\}-\{s_{i-1}\}  \in \bar{S}_k(U_{i-1})$$
is a cycle. 
Since $H_k(U_{i-1};\Z)\!=\!0$, it must be a boundary.
If $i\!=\!1$, this conclusion is still true 
with $U_0\!=\!X$, $K_0^{\top}\!=\!\eset$, and $s_0\!=\!0$, 
since $[h]\!=\!0$ by assumption. 
Let
$$\ti{s}_i\equiv\sum_{j=1}^{\ti{N}_i}\ti{f}_{i,j}\in S_{k+1}(U_{i-1})$$
be such that 
$$\sum_{\si\in K_i^{\top}-K_{i-1}^{\top}}\!\!\!\!\!\!
\{h\circ\eta\circ l_{\si}\}+\{s_i\}-\{s_{i-1}\} 
=\bar\partial\big\{\ti{s}_i\big\}\in \bar{S}_k(U_{i-1}).$$
Similarly to Lemma~\ref{gluing_l2}, there exist 
$$\ti\C_i^{(0)}\subset
\ti\C_i\!\equiv\!\bigsqcup_{i'=1}^{i'=i}\C_{\ti{s}_{i'}},$$
a symmetric subset $\ti\D_i\!\subset\!\ti\C_i\!\times\!\ti\C_i$ disjoint from 
the diagonal, and maps
\begin{gather*}
\ti\tau_i\!: \ti\D_i\lra {\cal S}_k, \quad 
\big((j_1,p_1),(j_2,p_2)\big)\lra\ti\tau_{i,((j_1,p_1),(j_2,p_2))},\\
(\ti{j}_i,\ti{p}_i)\!:K_i^{\top}\!\sqcup\!\{1,\ldots,N_i\}\lra\ti{\cal C}_i^{(0)},
~~~\hbox{and}~~~
\ti\tau_i\!:K_i^{\top}\!\sqcup\!\{1,\ldots,N_i\}
\lra{\cal S}_k, ~~j\lra\ti\tau_{(i,j)}, 
\end{gather*}
such that\\
${}\quad~$ (i) $\ti\D_i\!\subset\!\ti\D_{i+1}$, $\ti\tau_{i+1}|_{\ti\D_i}\!=\!\ti\tau_i$,
and $(\ti{j}_{i+1},\ti{p}_{i+1},\ti\tau_{i+1})|_{K_i^{\top}}=
(\ti{j}_i,\ti{p}_i,\ti\tau_i)|_{K_i^{\top}};$\\
${}~~~~$ (ii) the projection $\ti\D_i\!\lra\!\ti\C_i$  on either coordinate
is a bijection onto the complement of~$\ti\C_i^{(0)}$;\\
${}~~~$ (iii) for all $((j_1,p_1),(j_2,p_2))\!\in\!\ti\D_i\!\cap\!
(\C_{\ti{s}_{i_1}}\!\times\!\C_{\ti{s}_{i_2}})$,
\begin{gather*}
\ti\tau_{i,((j_2,p_2),(j_1,p_1))}=\ti\tau_{i,((j_1,p_1),(j_2,p_2))}^{~-1}, \qquad
\ti{f}_{i_2,j_2}\circ\io_{k+1,p_2}=
\ti{f}_{i_1,j_1}\circ\io_{k+1,p_1}\circ\ti\tau_{i,((j_1,p_1),(j_2,p_2))},\\
\hbox{and}\qquad \sign\ti\tau_{i,((j_1,p_1),(j_2,p_2))}=-(-1)^{p_1+p_2};
\end{gather*}
${}~~~$  (iv) for all $\si\!\in\!K_i^{\top}\!-\!K_{i-1}^{\top}$, 
\begin{equation*}
\ti{f}_{i,\ti{j}_i(j)}\circ\io_{k+1,\ti{p}_i(j)}\circ\ti\tau_{(i,j)}=f_{\si}
\qquad\hbox{and}\qquad \sign\ti\tau_{(i,j)}=-(-1)^{\ti{p}_i(j)};
\end{equation*}
${}~~~~$  (v) $(\ti{j}_i,\ti{p}_i)$ is a bijection onto $\ti\C_i^{(0)}$.\\

\noindent
(4) Put 
\begin{gather*}
\ti{M}'=\Big(\bigsqcup_{i=1}^{\i}\bigsqcup_{j=1}^{\ti{N}_i}
\{i\}\!\times\!\{j\}\!\times\!\De^{k+1}\sqcup I\!\times\!M'\Big)\Big/\sim, 
\qquad\hbox{where}\\
\big(i_1,j_1,\io_{k,p_1}(\ti\tau_{i,((j_1,p_1),(j_2,p_2))}(t))\big)\sim
 \big(i_2,j_2,\io_{k,p_2}(t)\big) \quad\forall~
 ((j_1,p_1),(j_2,p_2))\!\in\!\ti\D_i\!\cap\!
(\C_{\ti{s}_{i_1}}\!\times\!\C_{\ti{s}_{i_2}}),~t\!\in\!\De^k,\\
\big(1,\pi(\si,t)\big)\sim\big(i,\ti{j}_i(\si),
\io_{k+1,\ti{p}_i(\si)}(\ti\tau_{i,\si}(t))\big)
\quad\forall~t\!\in\!\De^k,~\si\!\in\!K_i^{\top}\!-\!K_{i-1}^{\top},~i\!\in\!\Z^+.
\end{gather*}
Let 
$$\ti\pi\!: \bigsqcup_{i=1}^{\i}\bigsqcup_{j=1}^{\ti{N}_i}
\{i\}\!\times\!\{j\}\!\times\!\De^{k+1}\sqcup I\!\times\!M'\lra \ti{M}'$$
be the quotient map.
Define 
\begin{equation*}
\ti{F}\!:\ti{M}'\lra X  \qquad\hbox{by}\qquad
\begin{aligned}
&\ti{F}\big([i,j,t]\big)=\ti{f}_{i,j}\big(\ti\vp_{k+1}(\vp_{k+1}(t))\big)
&\quad&\forall~t\!\in\!\De^{k+1},~j\!=\!1,\ldots,\ti{N}_i,~i\!\in\!\Z^+;\\
&\ti{F}[s,x]=F(x)&\quad&\forall~s\!\in\!I,~x\!\in\!M',
\end{aligned}
\end{equation*}
where $\ti\vp_{k+1}$ and $\vp_{k+1}$ are the self-maps of 
$\De^{k+1}$ provided by Lemma~\ref{bdpush_lmm}.
Similarly to the proof of Lemma~\ref{equiv1},
this map is well-defined and continuous.
Since the image~of 
$$\bigsqcup_{i=2}^{\i}\bigsqcup_{j=1}^{\ti{N}_i}
\{i\}\!\times\!\{j\}\!\times\!\De^{k+1} \sqcup
I\!\times\! \pi\Big(\bigsqcup_{i=2}^{\i}
\bigsqcup_{\si\in K_2^{\top}}\!\!\!\!\{\si\}\!\times\!\De^{k-1}\Big)$$
under $\ti{F}\!\circ\!\ti\pi$ is contained in the pre-compact subset $U_1$ of~$X$,
$\ti{F}(\ti{M}')$ is a pre-compact subset of~$X$ as well.\\

\noindent
Let $\ti{M}$ be the complement in~$\ti{M}'$ of the~set
$$\ti\pi\Big(\bigsqcup_{i=1}^{\i}\bigsqcup_{j=1}^{j=\ti{N}_i}
\{i\}\!\times\!\{j\}\!\times\!\ti{Y}\sqcup I\!\times\!(M\!-\!M')\Big),$$
where $\ti{Y}\!\subset\!\De^{k+1}$ is the $(k\!-\!1)$-skeleton.
Similarly to the proof of Lemma~\ref{equiv1}, 
$\Bd\ti{F}|_{\ti{M}}$ is of dimension at most $k\!-\!1$, 
$\ti{M}$ is a smooth oriented manifold boundary, 
$$\partial\ti{M}=-M,$$ 
$\ti{F}|_{\ti{M}}$ is a smooth map, and $\ti{F}|_M\!=\!F|_M$.
We note that in this case,
\begin{gather*} 
\ti{M}_0=I\!\times\!M-\{1\}\!\times\!\!\!\bigcup_{\si\in K_i^{\top}}
\!\!\!\!\pi\big(\{\si\}\!\times\!(\De^k\!-\!\Int\De^k)\big),\quad
\partial\ti{M}_0\approx -\{0\}\!\times\!M+
\{1\}\!\times\!\!\bigsqcup_{\si\in K^{\top}}\!\!\!\!\{\si\}\!\times\!\Int\!\De^k;\\
\ti{M}^*=\ti\pi\Big(\bigsqcup_{i=1}^{\i}\bigsqcup_{j=1}^{j=\ti{N}_i}
\{i\}\!\times\!\{j\}\!\times\!(\De^{k+1}\!-\!\ti{Y})\Big),
\qquad 
\partial\ti{M}^*=\ti\pi\Big(\bigsqcup_{i=1}^{\i}\bigsqcup_{\si\in K_i^{\top}}\!\!
\{i\}\!\times\!\{\ti{j}_i(\si)\}\!\times\!(\Int\De^{k+1}_{\ti{p}_i(\si)}\Big).
\end{gather*}
The topological space $\ti{M}$ is obtained from $\ti{M}^*$ and $\ti{M}_0$
by identifying 
$$\{1\}\!\times\!\!\bigsqcup_{\si\in K^{\top}}\!\!\!\!\{\si\}\!\times\!\Int\!\De^k$$
with  $\partial\ti{M}^*$.
We conclude that 
$$\partial\big(\ti{F}|_{\ti{M}}\big)=-F|_M,$$
i.e.~$F|_M$ and $h$ represent the zero element in ${\cal H}_k(M)$.\\

\vspace{.2in}

\noindent
{\it Department of Mathematics, SUNY, Stony Brook, NY 11790-3651\\
azinger@math.sunysb.edu}\\

\end{document}